\documentclass[11pt]{article}
\usepackage{amsmath}
\usepackage{amssymb}
\usepackage{amsthm,amsxtra}
\usepackage{epsf}
\usepackage{eepic}
\newlength{\mytopmargin}
\newlength{\myleftmargin}
\newtheorem{lemma}{Lemma}
\newtheorem{prop}[lemma]{Proposition}
\newtheorem{thm}{Theorem}
\newtheorem{cor}[thm]{Corollary}

\newcommand\psymmU{%
\begin{picture}(1,1)(0,0)%
\allinethickness{0.5pt}%
\path(0,0)(0,1)(1,1)(1,0)(0,0)%
\end{picture}}
\newcommand\psymmUU{%
\begin{picture}(1,1)(0,0)%
\allinethickness{0.5pt}%
\path(0,0)(0,1)(1,1)(1,0)(0,0)%
\put(0.5,0.5){\makebox(0,0){$\cdot$}}%
\end{picture}}
\newcommand\psymmO{%
\begin{picture}(1,1)(0,0)%
\allinethickness{0.5pt}%
\path(0,0)(0,1)(1,1)(1,0)(0,0)%
\path(0,0)(1,1)%
\end{picture}}
\newcommand\psymmS{%
\begin{picture}(1,1)(0,0)%
\allinethickness{0.5pt}%
\path(0,0)(0,1)(1,1)(1,0)(0,0)%
\path(1,0)(0,1)%
\end{picture}}
\newcommand\psymmu{%
\begin{picture}(1,1)(0,0)%
\allinethickness{0.5pt}%
\path(0,0)(0,1)(1,1)(1,0)(0,0)%
\path(0,0)(1,1)%
\path(0,1)(1,0)%
\end{picture}}

\newbox\tsymmUbox
\newbox\tsymmUUbox
\newbox\tsymmObox
\newbox\tsymmSbox
\newbox\tsymmubox
\setbox\tsymmUbox =\hbox{\kern0.75pt\setlength{\unitlength}{6pt}\psymmU \kern0.75pt}

\setbox\tsymmUUbox=\hbox{\kern0.75pt\setlength{\unitlength}{6pt}\psymmUU\kern0.75pt}
\setbox\tsymmObox =\hbox{\kern0.75pt\setlength{\unitlength}{6pt}\psymmO \kern0.75pt}
\setbox\tsymmSbox =\hbox{\kern0.75pt\setlength{\unitlength}{6pt}\psymmS \kern0.75pt}

\newbox\symmUbox
\newbox\symmUUbox
\newbox\symmObox
\newbox\symmSbox
\newbox\symmubox
\setbox\symmUbox =\hbox{\kern0.75pt\setlength{\unitlength}{4.5pt}\psymmU \kern0.75pt}
\setbox\symmUUbox=\hbox{\kern0.75pt\setlength{\unitlength}{4.5pt}\psymmUU\kern0.75pt}
\setbox\symmObox =\hbox{\kern0.75pt\setlength{\unitlength}{4.5pt}\psymmO \kern0.75pt}
\setbox\symmSbox =\hbox{\kern0.75pt\setlength{\unitlength}{4.5pt}\psymmS \kern0.75pt}
\setbox\symmubox =\hbox{\kern0.75pt\setlength{\unitlength}{4.5pt}\psymmu \kern0.75pt}

\def\symmO{{\copy\symmObox}}

\begin{document}

\title{\Large\bf Correlation Functions for Random Involutions}
\author{Peter J. Forrester$^{\dagger}$, Taro Nagao${}^*$ and 
Eric M. Rains$^{\ddagger}$}
\date{}
\maketitle

\begin{center}
\it
 $^{\dagger}$Department of Mathematics and Statistics, University of Melbourne, \\
Victoria 3010, Australia
\end{center}
\begin{center}
${}^*$
\it
Graduate School of Mathematics, Nagoya University, \\
Chikusa-ku, Nagoya 464-8602, Japan
\end{center} 
\begin{center}
\it
 $^{\ddagger}$ Department of Mathematics, University of California, 
Davis, CA 95616, USA
\end{center}

\bigskip
\begin{center}
\bf Abstract 
\end{center}
\par
\bigskip
\noindent 
Our interest is in the scaled joint distribution
associated with $k$-increasing subsequences
for random involutions with a prescribed number of fixed points.
We proceed by specifying in terms of correlation functions
the same distribution for a Poissonized model in which both the number
of symbols in the involution, and the number of fixed points, are
random variables.
From this, a de-Poissonization argument yields the scaled
correlations 
and distribution function
for the random involutions. These are found to coincide with the
same quantities
known in random matrix theory from the study of ensembles
interpolating between the orthogonal and symplectic universality
classes at the soft edge, the interpolation being due to a rank 1
perturbation.

\newpage
\section{Introduction}
To motivate our study of random involutions, we first recall the
corresponding problem for random permutations, the solution of
which is known. 
Let
$S_N$ denote the set of the $N!$ distinct permutations of 
$\{1,2,\dots,N\}$. For each $\pi \in S_N$ denote the image of the number
$i \in \{1,2,\dots,N\}$ by $\pi(i)$. A subsequence of image points
$\pi(i_1), \pi(i_2), \dots, \pi(i_j)$ where $1 \le i_1 < \cdots < i_j \le N$
is said to be an increasing subsequence of length $j$ if
$\pi(i_1) < \pi(i_2) < \cdots < \pi(i_j)$. More generally, we say there is a
$k$-increasing subsequence of length $j$ if $\pi$ contains $k$ 
disjoint subsequences
$$
\pi(i_1^{(l)}) < \pi(i_2^{(l)}) < \cdots < \pi(i_{j_l}^{(l)}), \quad
1 \le i_1^{(l)} < \cdots < i_{j_l}^{(l)} \le N \: \: (l=1,\dots,k)
$$
with $\sum_{l=1}^k j_l = j$. Being disjoint these subsequences 
contain no common member. For a given $\pi$, let $L_N^{(k)}(\pi)$ denote the
maximum length of all the $k$-increasing subsequences and define
\begin{equation}\label{1.b}
\lambda_N^{(k)}(\pi) = L_N^{(k)}(\pi) - L_N^{(k-1)}(\pi), \qquad
L_N^{(-1)}(\pi) := 0.
\end{equation}
Note that
$$
\lambda_N^{(1)}(\pi) \ge \lambda_N^{(2)}(\pi) \ge \cdots \ge 
\lambda_N^{(N)}(\pi) \ge 0.
$$
Consider an ensemble of 
permutations of $\{1,2,\dots,N\}$ in which each permutation is equally
likely. The problem of interest is the computation of the scaled
joint distribution of $\{\lambda_N^{(j)}(\pi)\}_{j=1,\dots,l}$ in the limit
$N \to \infty$.

In the case $l=1$ it was proved by Baik, Deift and Johansson 
\cite{BDJ99} that
\begin{equation}\label{n2}
\lim_{N \to \infty} {\rm Pr} \Big ( {\lambda_N^{(1)} - 2 \sqrt{N} \over
N^{1/6} } \le s \Big ) = F_2(s),
\end{equation}
where $F_2(s)$ is the scaled cumulative distribution of the largest
eigenvalue for large random Hermitian matrices with complex elements
(technically matrices from the Gaussian unitary ensemble (GUE))
\cite{Fo93a,TW94a}. 
The latter is specified in terms of a Fredholm determinant according to
$$
F_2(s) = \det (1 - K^{\rm soft} \chi_{(s,\infty)}),
$$
where $\chi_J$ is the characteristic function of the interval $J$, while
$K^{\rm soft}$ is the integral operator with kernel given in terms of
Airy functions according to
\begin{eqnarray}\label{1.2so}
 K^{\rm soft}(x,y) & = & {{\rm Ai}(x) {\rm Ai}'(y) - {\rm Ai}(y)
{\rm Ai}'(x) \over x - y} \nonumber \\
& = & \int_0^\infty {\rm Ai}(x+t) {\rm Ai}(y+t) \, dt.
\end{eqnarray}
(The superscript `soft' indicates that the eigenvalue density
$\rho_{(1)}(s)$ is not strictly zero for any $s$, even though the origin
is chosen in the neighbourhood of the largest eigenvalue and thus at the
spectrum edge). The case of general $l$ was solved in \cite{BOO,KJ2,O00}, 
where it was proved that 
$$
\lim_{N \rightarrow \infty} 
{\rm Pr} \Big ( {\lambda_N^{(1)} - 2 \sqrt{N} \over N^{1/6}} \le s_1,
\dots, {\lambda_N^{(l)} - 2 \sqrt{N} \over N^{1/6}} \le s_l \Big )
= F_2(s_1,\dots,s_l).
$$
Here $F_2(s_1,\dots,s_l)$ is the scaled joint distribution of the $l$
largest eigenvalues in the GUE. The latter is uniquely specified in
terms of the $k$-point correlation function
$$
\rho_{(k)}^{\rm scaled}(s_1,\dots,s_k) = 
\det [ K^{\rm soft}(s_j,s_l) ]_{j,l=1,\dots,k}
$$
for the scaled eigenvalues at the soft edge of the GUE.

It is the objective of this study to calculate the analogous joint
probability 
for involutions with a prescribed number of fixed points.
One recalls that an involution is a permutation $\pi$ with the additional
property that $\pi^2 = I$, where $I$ denotes the identity permutation.
Involutions must consist entirely of two cycles and fixed points. Thus for
an involution of $\{1,2,\dots,N\}$, if there are $n$ two cycles, there must
be $m=N-2n$ fixed points. As emphasized in \cite{BR02a}, a random involution
with a prescribed number of fixed points can be generated geometrically by
marking $n$ points $(n \le [N/2])$
in the unit square below the diagonal $y=x$ uniformly at
random, marking the images of these points under reflection about $y=x$,
and marking $N-2n$ points uniformly at random on the diagonal. 
Equivalently the unit square can first be divided into a $N \times N$
integer grid, and the points marked at random on the lattice sites
below and on the diagonal according to the above prescription, with the
additional constraint that no two points are in the same row or column.
Either way,
projecting
the points onto the $x$-axis gives a sequence of $x$ co-ordinates
$x_1<x_2<\cdots<x_N$ while projecting them onto the $y$-axis gives a
sequence of $y$ co-ordinates $y_1 < y_2 < \cdots < y_N$. Each point
will then have a co-ordinate $(x_i,y_{\pi(i)})$ with the property that
$\pi(i) = i$ for the point on the diagonal, and that $(x_{\pi(i)},y_i)$
is the point reflected in the diagonal otherwise. Hence $\pi$ defines an
involution with a prescribed number of fixed points, and furthermore
the involutions are generated at random with uniform probability by this
procedure.

\begin{figure}[t]
\epsfxsize=7cm
\centerline{\epsfbox{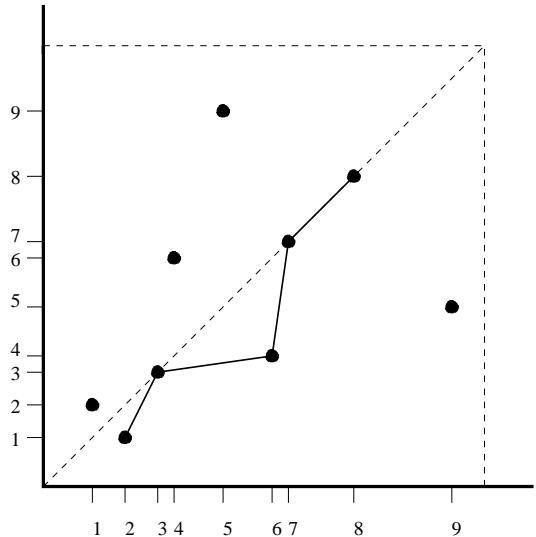}}
\caption{\label{f.inv.1} Nine points in the square symmetric about the
diagonal which correspond to the involution $(12)(3)(46)(59)(7)(8)$.
An rd path realizing $L_N^{(1)}(\pi)$ has been drawn. } 
\end{figure}

The quantity $L_N^{(k)}(\pi)$ admits an interpretation in the above setting.
Connect points by segments which always have positive slope to form a
continuous path, which is said to be right/diagonal (rd).
Define the length of this path
as the number of points it contains, and denote it by $\#$rd. Similarly,
let $({\rm RD})^k$ denote the set  of all $k$ disjoint rd lattice paths,
formed from amongst the points with a path being disjoint if it contains
no common points, and for $({\rm rd})^k \in ({\rm RD})^k$ let
$\# ({\rm rd})^k$ denote the number of lattice points. Then by considering
recurrences satisfied by the various quantities one can show
\cite{Ra02}
\begin{equation}\label{n4}
L_N^{(k)}(\pi) = {\rm max} \sum_{({\rm rd})^k \in ({\rm RD})^k}
\# ({\rm rd})^k
\end{equation}
(for the recurrences satisfied by the left hand side, see \cite[Appendix A]{FR02b}).
The equation (\ref{n4}) holds for points corresponding to a general
permutation $\pi$. In the case that $\pi$ is an involution, it is easy
to see that the set of lattice paths $({\rm RD})^k$ can be restricted to those
which contain points on or below the diagonal (see Figure \ref{f.inv.1} for an
example).

The explicit form of the scaled distribution of 
$L_N^{(1)} =: \lambda_{n,m}^{(1)}$ for
random involutions with $n$ 2-cycles and $m$ fixed points,
and thus the solution to our problem in the case $l=1$,
is already known \cite{BR02b}. Thus
introduce the scaling variable $w$ by the requirement that
\begin{equation}\label{s1}
m = [ \sqrt{2n} - 2w (2n)^{1/3} ],
\end{equation}
where $[ \cdot ]$ denotes the integer part.
It is proved in \cite{BR02b}
that with $w$ fixed
\begin{equation}\label{s2}
\lim_{N \to \infty} {\rm Pr} \Big ( {\lambda_{n,m}^{(1)} - 2 \sqrt{N} \over
N^{1/6} } \le s \Big ) = F^{\symmO}(s;w),
\end{equation}
where the distribution $F^{\symmO}(s;w)$ is specified in terms of a certain
Riemann-Hilbert problem related to the Painlev\'e II  equation with special
monodromy data. This distribution has the property that
$$
F^{\symmO}(s;0) = F_1(s), \qquad \lim_{w \to \infty} 
F^{\symmO}(s;w) = F_4(s),
$$
where $F_1(s), F_4(s)$ are the scaled cumulative distributions for the
largest eigenvalue in the GOE of random real symmetic matrices, and the
GSE of Hermitian matrices with real quaternion elements respectively.

Here we seek the joint distribution of
$\{\lambda_{n,m}^{(j)}(\pi)\}_{j=1,\dots,l}$, scaled as in (\ref{s1}),
(\ref{s2}) for general $l$. Regarding the joint distribution of 
all the $\lambda_{n,m}^{(j)}(\pi)$, $j=1,\dots,N$
as specifying a point process, it is
generally true that the 
distribution of the $l$ right-most points is fully determined
by the corresponding correlation functions. 
The specification of the correlations for a point process
implying the distribution  $F^{\symmO}(s;w)$ for the right-most
point has recently been given \cite{FR02a}. This was found
in the study of 
a closely related geometrical model to that corresponding to random 
involutions and their increasing subsequences.
Thus consider an $M \times M$ integer grid. Associate
with each lattice site on or below the diagonal a continuous
exponential variable
\begin{eqnarray}\label{2.1g}
{\rm Pr}(x_{i,j} \in [y,y+dy]) & = & e^{-y} \, dy, \qquad i < j, \nonumber \\
{\rm Pr}(x_{i,i} \in [y,y+dy]) & = & {(1-A) \over 2}
e^{-(1-A)y/2} \, dy,
\end{eqnarray}
and impose the symmetry constraint that $x_{i,j} = x_{j,i}$ for $i > j$.
The quantities (\ref{n4}), with the lattice points counted according to their weightings
$x_{i,j}$, are well defined. Setting then
$x_k = L_M^{(k)} - L_M^{(k-1)}$, $L_M^{(-1)} :=0$ we have
$x_1 > x_2 > \cdots > x_M > 0$ and furthermore the joint distribution of
these variables is proportional to \cite[Section 3]{Ba03} \cite[Prop.~4]{FR02b}
\begin{equation}\label{xj}
e^{-\sum_{j=1}^M x_j/2}
e^{A \sum_{j=1}^M (-1)^{j-1} x_j/2} \prod_{1 \le i < j \le M}
(x_i - x_j).
\end{equation}

This same p.d.f.~occurs in random matrix theory. Thus let $X$ be a
$2n \times 2n$ antisymmetric complex Gaussian matrix
(independent entries distributed
as ${\rm N}[0,1]+i {\rm N}[0,1]$), and let $\vec{x}$ be
a $2n \times 1$ complex Gaussian vector with entries 
distributed as
${\rm N}[0,1/\sqrt{2}]+i {\rm N}[0,1/\sqrt{2}]$. Then we know from
\cite[Thm.~3]{FR02b} that $Y = X^\dagger X + b \vec{x} \vec{x}^\dagger$
has eigenvalue p.d.f.~(\ref{xj}) with $A = 1-2/b$, $M=2n$.

For the p.d.f.~(\ref{xj}), the scaled correlation functions with
\begin{equation}\label{Sc}
x_j = 4M + 2(2M)^{1/3} X_j, \qquad A= u/2(2M)^{1/3},
\end{equation}
and $M \to \infty$ have been computed \cite{FR02a} (actually the symbol
$\alpha$ is used instead of $u$ in \cite{FR02a}, but we use $\alpha$
for another purpose below).  On the other
hand, a result of Baik \cite{Ba03} gives
\begin{equation}\label{Ba1}
\lim_{M \to \infty}
{\rm Pr} \Big ( {L_M^{(1)} |_{A = -2w/(2M)^{1/3}} - 4M \over 2 (2M)^{1/3}}
\le s \Big ) = F^{\symmO}(s;w).
\end{equation} 
The scalings (\ref{Sc}) is precisely that in (\ref{Ba1}) with the
identification
\begin{equation}\label{1.10a}
w = - {u \over 4},
\end{equation}
so this
distribution is fully determined by the correlations computed in
 \cite{FR02a}. We therefore expect that for random involutions, the
generalization of (\ref{s2}) is 
\begin{equation}\label{ff}
\lim_{N \to \infty} {\rm Pr} \Big ( {\lambda_{n,m}^{(1)} - 2 \sqrt{N} \over
N^{1/6} } \le s_1, \dots,
{\lambda_{n,m}^{(l)} - 2 \sqrt{N} \over
N^{1/6} } \le s_l
 \Big ) = F^{\symmO}(s_1,\dots,s_l;w), 
\end{equation}
where $F^{\symmO}(s_1,\dots,s_l;w)$ is the joint distribution of the
scaled $l$ right-most points in the process specified by
(\ref{xj}), scaled as in (\ref{Sc}) and with $u$ and $w$ related by
(\ref{1.10a}). As revised below, it is generally true that the latter
distribution is fully determined by the scaled correlation functions
$\rho_{(k)}^{\rm scaled}$. According to \cite{FR02a} these have
the explicit form
\begin{equation}\label{1.11}
\rho_{(k)}^{\rm scaled}(X_1,\dots,X_k;u) = {\rm qdet} \,
[ f(X_i,X_j) ]_{i,j=1,\dots,k},
\end{equation}
where with ${\rm sgn}(u)=1$, $-1$, 0 for $u>0$, $u<0$, $u=0$ respectively,
$f$ is a $2 \times 2$ matrix with entries
\begin{eqnarray}\label{1.12}
f^{11}(X,Y) & = & f^{22}(Y,X), \nonumber \\
f^{22}(X,Y) & = & {1 \over 2} K^{\rm soft}(X,Y) - {1 \over 2}
{\partial \over \partial Y} \int_{-\infty}^X e^{u(X-t)/2}
K^{\rm soft}(t,Y) \, dt \nonumber \\
&& - {u \over 4} \int_{-\infty}^X dt \, e^{u(X-t)/2}
\int_Y^\infty ds \, {\partial \over \partial t} K^{\rm soft}(s,t),
 \nonumber \\
f^{12}(X,Y) & = & {1 \over 4} 
\Big ( {u \over 2} + {\partial \over \partial X} \Big )
\Big ( {u \over 2} + {\partial \over \partial Y} \Big )
\Big \{ \int_X^{\infty} K^{\rm soft}(Y,t) \, dt -
 \int_Y^{\infty} K^{\rm soft}(X,t) \, dt \Big \}, \nonumber \\
f^{21}(X,Y) & = & - e^{u|X-Y|/2} {\rm sgn}(X-Y) - \Big \{
\int_{-\infty}^Y e^{u(Y-t)/2} K^{\rm soft}(X,t) \, dt  \nonumber \\
&& - \int_{-\infty}^X e^{u(X-t)/2}  K^{\rm soft}(Y,t) \, dt 
\Big \},
\end{eqnarray}
and the qdet operation is defined in terms of the more familiar
Pfaffian by the general formula
$$
{\rm qdet} \, A = {\rm Pf} (A Z_{2n}^{-1}), \qquad
Z_{2n} = 1_n \otimes
\left [ \begin{array}{cc} 0 & -1 \\ 1 & 0 \end{array} \right ],
$$
valid for all $2n \times 2n$ matrices A with the self dual property  
$A = Z_{2n}^{-1} A^T Z_{2n}$. Our main task then is to show that the
same scaled correlations determine the joint distribution of
$\{\lambda_{n,m}^{(j)}(\pi)\}_{j=1,\dots,l}$. 
We will see that
doing this indeed allows (\ref{ff}) to be validated,
giving us the following 
limit theorem, which is our main result.

\begin{thm}
Consider an involution on $N = 2n+m$ symbols, consisting of $m$ fixed
points and $n$ 2-cycles, chosen at random with uniform distribution on
the set of such involutions. Let $\lambda_{n,m}^{(k)}$ be specified in
terms of the maximum length of the $k$-increasing subsequences
$L_N^{(k)}$ according to (\ref{1.b}). With $m$ related to $n$ by
(\ref{s1}), the limit formula (\ref{ff}) holds.
\end{thm}

\section{Random Generalized Involutions}
\setcounter{equation}{0}

\subsection{Strategy}
We are guided in our approach by previous studies on the calculation of
distribution functions relating to increasing subsequences for random
permutations \cite{BDJ99,BOO,KJ1}, and also previous studies on the
calculation of the distribution of the maximum increasing subsequence for
random involutions \cite{BR02b,BF03}. For both random permutations and
random involutions, all the studies proceed by Poissonizing the ensemble.
In the case of random permutations, the number of symbols $N$ is
itself taken as a variable which occurs with probability
$e^{-z} z^N/ N!$. For random involutions the number of two cycles $n$
and number of fixed points $m$ are separately taken as random variables, 
occurring with probability $e^{-z_1 - z_2} z_1^n z_2^m / n! m!$. 
One then seeks to
calculate the distribution function of interest in the Poissonized
ensemble, and to compute the $z \to \infty$ ($z_1, z_2 \to
\infty$) scaled limit. A
de-Poissonization argument \cite{Jo98} gives that the latter is
equivalent to the $N \to \infty$ ($n,m \to \infty$)
limit in the original ensemble. 

As the distribution functions are controlled by the correlation functions
$\rho_{(k)}$,
one would like to first calculate the Poissonized correlation functions.
However we know from studies of random permutations that a direct
calculation is not practical. Instead, as was shown by Johansson \cite{KJ1},
progress can be made by constructing the Poissonization as a limiting case of
a generalization of the geometrical viewpoint of the increasing
subsequence problem. In one such generalization, giving rise to the so
called Meixner ensemble \cite{KJ2}, each lattice point of an $M \times M$
grid carries a non-negative integer variable chosen from the geometric
distribution with parameter $q$. For this model $\rho_{(k)}$ can be
computed in terms of certain orthogonal polynomials. Taking the limit
$M \to \infty$ with $q=Q/M^2$ gives $\rho_{(k)}$ for the
Poissonized version of the original model. Thus we must first calculate
$\rho_{(k)}$ for the version of the  Meixner ensemble relevant to
involutions.

\subsection{The joint p.d.f., correlations and distribution
functions}\label{spdf}
According the above strategy, our first task is to introduce a
generalization of the geometrical model corresponding to random
involutions and their increasing subsequences. For this we consider
an $M \times M$ grid, and associate with each lattice site on or below
the diagonal a non-negative integer variable
\begin{eqnarray}\label{jp}
{\rm Pr} \, (x_{i,j} = k) & = & (1 - q) q^k, \qquad i < j, \nonumber \\
{\rm Pr} \, (x_{i,i} = k) & = & (1 - \sqrt{\alpha q}) (\alpha q)^{k/2}.
\end{eqnarray}
For $i > j$ we impose the symmetry constraint $x_{i,j} = x_{j,i}$.
Notice that the symmetry constraints of this model are the same as that
for the geometrical model of random involutions specified in the
Introduction, and notice too the similarity with the model defined by
(\ref{2.1g}) and surrounding text. As with the latter model, the quantities
(\ref{n4}), which can be regarded as a sequence of last passage times,
are well defined. With $\lambda_k := L_M^{(k)} - L_M^{(k-1)}$,
$L_M^{(-1)} := 0$, we know from 
\cite[Section 3]{Ba03}, \cite[Prop.~1]{FR02b} that the joint distribution of the latter
variables have the explicit form
\begin{equation}\label{1.9P}
P_M(\lambda;q) = (1 - \sqrt{\alpha q})^M (1 - q)^{M(M-1)/2}
 q^{\sum_{j=1}^M \lambda_j/2}
\alpha^{\sum_{j=1}^M(-1)^{j-1} \lambda_j/2}
\prod_{1 \le j < l \le M}
{\lambda_j - \lambda_l + l - j \over l - j}.
\end{equation}
In terms of 
$h_j := \lambda_j + M- j$ this reads
\begin{equation}\label{1.9}
P(h_1,\cdots,h_M) = C_M(q,\alpha)  q^{\sum_{j=1}^M h_j/2}
\alpha^{\sum_{j=1}^M (-1)^{j-1} h_j}
\prod_{1 \le j < l \le M}  (h_j - h_l),
\end{equation}
where
\begin{equation}\label{2.10}
C_M(q,\alpha) = (1 - \sqrt{\alpha q})^M (1 - q)^{M(M-1)/2}
q^{-M(M-1)/2}  \alpha^{-M/4} \prod_{j=1}^{M-1} (1/j!)
\end{equation}
and $\infty > h_1 > h_2 > \cdots > h_M \geq 0$. Note the similarity
between (\ref{1.9}) and (\ref{xj}).

Let us introduce the symmetrized joint p.d.f.~by
\begin{equation}\label{9.2}
P_{\rm sym}(h_1,\dots,h_M) = \sum_{\mu \in S_M}
P(h_{\mu(1)},\dots, h_{\mu(M)})
\chi_{h_{\mu(1)} > \cdots > h_{\mu(N)}},
\end{equation}
where $\chi_T = 1$ for $T$ true, and $\chi_T = 0$ otherwise. The $k$-point
correlation $\rho_{(k)}$ is then given by
\begin{equation}\label{2.5a}
\rho_{(k)}(h_1,\dots,h_k) =
{1 \over (M-k)!} \sum_{h_{k+1},\dots,h_M = 0}^\infty
P({h}_1,\dots,h_M).
\end{equation}

The correlations (\ref{2.5a}) can be used to compute the distribution
functions ${\rm Pr}(h_1 \le a_1,\dots,h_l \le a_l)$ for the joint
p.d.f.~(\ref{1.9}). To see this, let $a_0 = \infty$, let
$a_1 > a_2 > \cdots > a_l$ be non-zero integers, and put
$I_j = (a_j,a_{j-1})$ where $(a_j, a_{j-1})$ denotes all integers
between (but not including) $a_j$ and $a_{j-1}$. Let
$E_M(\{(n_r,I_r) \}_{r=1,\dots,l})$ denote the probability that
$n_r$ of the coordinates $\{h_i\}_{i=1,\dots,M}$ are in $I_r$
$(r=1,\dots,l)$. Then as a consequence of the definitions we know
from \cite[eq.~(3.41)]{KJ2} that with
$$
\mathbb L_l := \{(n_1,\dots,n_l) \in \mathbb Z_{\ge 0}^l \: : \:
\sum_{j=1}^r n_j \le  r-1 \: (r=1,\dots,l) \},
$$
one has
\begin{equation}\label{6a.c}
{\rm Pr}(h_1 \le a_1, \dots, h_l \le a_l) =
\sum_{(n_1,\dots,n_l) \in \mathbb L_l}
E_M(\{(n_r,I_r) \}_{r=1,\dots,l}).
\end{equation}
Furthermore, again from the definitions, it is easy to see that
\begin{equation}\label{6a.1}
E_M(\{(n_r,I_r) \}_{r=1,\dots,l})
= {(-1)^{\sum_{r=1}^l n_r} \over n_1! \cdots n_l !}
{\partial^{\sum_{j=1}^l n_j} \over
\partial \xi^{n_1} \cdots \partial \xi^{n_l} }
\Big \langle \prod_{j=1}^M \Big (1 - \sum_{r=1}^k \xi_r \chi_{I_r}^{(j)}
\Big ) \Big \rangle_{\rm Psym} \Big |_{\xi_1=\cdots=\xi_k=1},
\end{equation}
where $\chi_{I_r}^{(j)} = 1$ if $h_j \in I_r$,
$\chi_{I_r}^{(j)} = 0$ otherwise, and that the average herein is given
in terms of the correlations (\ref{2.5a}) according to
\begin{equation}\label{6a.2}
\Big \langle \prod_{j=1}^M \Big (1 - \sum_{r=1}^k \xi_r \chi_{I_r}^{(j)}
\Big ) \Big \rangle_{\rm Psym} =
1 + \sum_{p=1}^M {(-1)^p \over p!}
\sum_{h_1,\dots,h_p=0}^\infty
\prod_{i=1}^p \Big ( \sum_{r=1}^k \xi_r \chi_{I_r}^{(i)} \Big )
\rho_{(p)}(h_1,\dots,h_p).
\end{equation}

\subsection{Quaternion determinant expression for the $k$-point correlation}
We know from \cite[Eq.~(3.1)]{FR02a} that with $P$ given by (\ref{1.9}), and $M$ even
(\ref{9.2}) can be written 
\begin{equation}
P_{\rm sym}(h_1,\dots,h_M) = C_M(q,\alpha)q^{\sum_{j=1}^M h_j/2}
\prod_{j>l}^M (h_j - h_l) {\rm Pf}[\epsilon(h_j,h_l)]_{j,l=1,2,\cdots,M},
\end{equation}
where
\begin{equation}\label{ep}
\epsilon(x,y) = \alpha^{|y-x|/2} {\rm sgn}(y-x)
\end{equation}
(a similar formula can be given in the case that $M$ odd, but for
definiteness we will proceed with the assumption that $M$ is even).
This structure is familiar in random matrix theory, and in the
theory of random measures on partitions relating to increasing
subsequences \cite{RAINS}. It is known \cite{FP95,MAHOUX} that 
in general the corresponding $k$-point correlations have the 
quaternion determinant form
\begin{equation}
\rho_{(k)}(h_1,\cdots,h_k) =
\ {\rm qdet}[f(h_j,h_l)]_{j,l=1,2,\cdots,k}, 
\end{equation}
where $f(x,y)$ is the $2 \times 2$ matrix
\begin{equation}\label{fxy}
f(x,y) = \left[ \begin{array}{cc} S(x,y) & I(x,y) \\
D(x,y) & S(y,x) \end{array} \right].
\end{equation}
According to the general formalism, 
to specify the matrix elements in (\ref{fxy})
we must
introduce
skew orthogonal polynomials $R_n(x)$ satisfying
\begin{displaymath}
\langle R_{2m}(x), R_{2n+1}(y) \rangle =
- \langle R_{2n+1}(x), R_{2m}(y) \rangle = r_m \delta_{m,n},
\end{displaymath}
\begin{equation}\label{3.8}
\langle R_{2m}(x), R_{2n}(y) \rangle = 0, \quad
 \langle R_{2m+1}(x), R_{2n+1}(y) \rangle = 0,
\end{equation}
where
\begin{eqnarray}
 \langle f(x), g(y) \rangle & = &
\frac{1}{2} \sum_{x=0}^{\infty} \sum_{y=0}^{\infty}
q^{(x+y)/2} \epsilon(x,y) (f(x)g(y) - f(y)g(x))
\nonumber \\ &  = &
\sum_{y=0}^{\infty} q^{y/2} \alpha^{-y/2} f(y) \sum_{x=y+1}^\infty
q^{x/2} \alpha^{x/2} g(x)
- \sum_{y=0}^{\infty} q^{y/2} \alpha^{-y/2}
g(y) \sum_{x=y+1}^{\infty} q^{x/2} \alpha^{x/2} f(x).
\nonumber \\
\label{antib}
\end{eqnarray}
One then has
\begin{equation}\label{2.11}
f(x,y) = \left[ \begin{array}{cc} S(x,y) & I(x,y) \\
D(x,y) & S(y,x) \end{array} \right],
\end{equation}
where
\begin{equation}
S(x,y) = q^{y/2} \sum_{j=0}^{(M/2)-1} \frac{1}{r_j} \left[
\Phi_{2 j}(x) R_{2 j + 1}(y) - \Phi_{2 j + 1}(x) R_{2 j}(y) \right],
\label{sxy}
\end{equation}
\begin{equation}\label{ixy}
I(x,y) = - \sum_{j=0}^{(M/2)-1} \frac{1}{r_j} \left[
\Phi_{2 j}
(x) \Phi_{2 j + 1}(y) - \Phi_{2 j + 1}(x)
\Phi_{2 j}(y) \right]
+ \epsilon(x,y)
\end{equation}
and
\begin{equation}
D(x,y) = q^{(x+y)/2} \sum_{j=0}^{(M/2)-1} \frac{1}{r_j} \left[
R_{2 j}(x) R_{2 j + 1}(y) - R_{2 j + 1}(x) R_{2 j}(y) \right]
\label{dxy}
\end{equation}
with
\begin{equation}\label{phj}
\Phi_j(x) = \sum_{y=0}^{\infty} \epsilon(y,x) q^{y/2} R_j(y).
\end{equation}

\subsection{Skew orthogonal polynomials}
The explicit forms of polynomials with the skew orthogonality
property (\ref{3.8}) with respect to the skew
product (\ref{antib}) is required. These can be calculated from
general determinant formulas
for skew orthogonal polynomials.

\begin{lemma}
Let $\langle \, , \, \rangle$ be a general skew symmetric product. Let
$\{ R_j(x) \}_{j=0,1,\dots}$ be the corresponding skew orthogonal
polynomials which thus satisfy (\ref{3.8}). Let $\{ C_j(x) \}_{j=0,1,\dots}$
be any family of monic polynomials. Assuming ${\cal D}_n$ and  
${\cal E}_n$ as specified by (\ref{dn}) and (\ref{en}) below
are non-zero we have
\begin{equation}\label{R1}
R_{2 n}(x) = {\cal D}_n^{-1} \left| \begin{array}{ccccc}
C_{2 n}(x) & J^{2n \ 2n-1} & J^{2n \ 2n-2} & \cdots & J^{2n \ 0} \\
C_{2 n - 1}(x) & J^{2n-1 \ 2n-1} & J^{2n-1 \ 2n-2} & \cdots & J^{2n-1 \ 0} \\
\vdots & \vdots & \vdots & \ddots & \vdots \\
C_0(x) & J^{0 \ 2n-1} & J^{0 \ 2n-2} & \cdots & J^{0 \ 0} \end{array} \right|
\end{equation}
and
\begin{eqnarray}\label{R2}
R_{2 n+1}(x) & = & {\cal E}_n^{-1} \left| \begin{array}{cccccc}
J^{2n+1 \ 2n+1} & C_{2 n+1}(x) & J^{2n+1 \ 2n-1}
& J^{2n+1 \ 2n-2} & \cdots & J^{2n+1 \ 0} \\
J^{2n \ 2n+1} & C_{2 n}(x) & J^{2n \ 2n-1}
& J^{2n \ 2n-2} & \cdots & J^{2n \ 0} \\
\vdots & \vdots & \vdots & \vdots & \ddots & \vdots \\
J^{0 \ 2n+1} & C_0(x) & J^{0 \ 2n-1} & J^{0 \ 2n-2} &
\cdots & J^{0 \ 0} \end{array} \right|  \nonumber \\
&  & + c_n R_{2n}(x),
\end{eqnarray}
where
\begin{equation}\label{3.9}
J^{mn} = \langle C_m(x),C_n(y) \rangle,
\end{equation}
\begin{equation}\label{dn}
{\cal D}_n = \left| \begin{array}{cccc}
J^{2n-1 \ 2n-1} & J^{2n-1 \ 2n-2} & \cdots & J^{2n-1 \ 0} \\
J^{2n-2 \ 2n-1} & J^{2n-2 \ 2n-2} & \cdots & J^{2n-2 \ 0} \\
\vdots & \vdots & \ddots & \vdots \\
J^{0 \ 2n-1} & J^{0 \ 2n-2} & \cdots & J^{0 \ 0} \end{array} \right|,
\end{equation}
\begin{equation}\label{en}
{\cal E}_n = - \left| \begin{array}{cccc}
J^{2n \ 2n+1} & J^{2n \ 2n-1} & \cdots & J^{2n \ 0} \\
J^{2n-1 \ 2n+1} & J^{2n-1 \ 2n-1} & \cdots & J^{2n-1 \ 0} \\
\vdots & \vdots & \ddots & \vdots \\
J^{0 \ 2n+1} & J^{0 \ 2n-1} & \cdots & J^{0 \ 0} \end{array} \right|
\end{equation}
and the $c_n$ are arbitrary constants.
\end{lemma}

\noindent
Proof. \quad We note the skew orthogonality properties (\ref{3.8}) are
equivalent to the requirements
\begin{equation}\label{11.1}
\langle C_j(x), R_{2n}(y) \rangle = \langle C_j(x), R_{2n+1}(y) \rangle
= 0, \qquad j \le 2n - 1
\end{equation}
and
\begin{equation}\label{2.22a}
\langle C_{2n}(x), R_{2n+1}(y) \rangle \ne 0.
\end{equation}
 The property
(\ref{11.1}) is immediate from the structure of (\ref{R1}) and
(\ref{R2}), while (\ref{2.22a}) follows from the assumption
that ${\cal D}_{n+1}$ is non-zero. Furthermore, by inspection both
(\ref{R1}) and (\ref{R2}) give monic polynomials.
\hfill $\square$

\medskip
The structure of (\ref{antib}) suggests choosing $\{C_j(x)\}_{j=0,1,\dots}$
as particular Meixner polynomials. We recall \cite{KOEKOEK} the
Meixner polynomials are defined as
\begin{equation}\label{3.1}
M_n(x;c,q)= \ _2F_1 \left( \begin{array}{c} -n,-x \\ c \end{array}
\right| \left. 1 - \frac{1}{q} \right)
\end{equation}
and satisfy
\begin{equation}\label{2.24}
\sum_{x=0}^{\infty} \frac{(c)_x}{x!} q^x M_m(x;c,q) M_n(x;c,q)
= \frac{q^{-n} n!}{(c)_n (1 - q)^c} \delta_{m,n}
\end{equation}
with $(c)_n = \Gamma(c + n)/\Gamma(c)$.
We define monic
polynomials $C_n(x)$ as
\begin{equation}\label{3.3}
C_n(x) = \frac{n! q^n}{(q-1)^n} M_n(x;1,q) = x^n + \cdots.
\end{equation}
From (\ref{2.24}) the orthogonality relation for the $C_n(x)$ is 
\begin{equation}\label{ortho1}
\sum_{x=0}^{\infty} q^x C_m(x) C_n(x) = h_n \delta_{m,n},
\end{equation}
where
\begin{equation}\label{3.6}
h_n = \frac{(n!)^2 q^n}{(1 - q)^{2 n + 1}}.
\end{equation}
Using the symmetry $M_n(x;\gamma,q) = M_x(n;\gamma,q)$,
(\ref{2.24}) also implies the completeness relation
\begin{equation}
\sum_{n=0}^{\infty} \frac{1}{h_n} C_n(x) C_n(y) =
q^{-x} \delta_{x,y}, \qquad x,y \in  \mathbb Z_{\ge 0}.
\label{ortho2}
\end{equation}

With $\{C_j(x) \}$ so specified we seek to compute (\ref{3.9}).

\begin{prop}
For $m < n$,
\begin{equation}\label{3.15}
J^{mn} = - J^{nm} = \langle C_m(x),C_n(y) \rangle  =
a_m b_n,
\end{equation}
where
\begin{eqnarray}\label{3.15a}
a_n & = & n! \frac{\sqrt{\alpha}}{\sqrt{\alpha} - \sqrt{q}}
\left(\frac{\sqrt{q} ( 1 - \sqrt{\alpha q})}{(1-q)(\sqrt{\alpha} -
\sqrt{q})} \right)^n, \nonumber \\
b_n & = & n! \frac{1}{1 - \sqrt{\alpha q}}
\left(\frac{\sqrt{q} ( \sqrt{\alpha} - \sqrt{q})}{(1-q)(1 -
\sqrt{\alpha q})} \right)^n.
\end{eqnarray}
\end{prop}

\noindent
Proof. \quad Set
\begin{equation}\label{2.30a}
F(y;\alpha) =
q^{-y/2} \alpha^{-y/2} \sum_{x=y+1}^\infty q^{x/2} \alpha^{x/2}
C_n(x).
\end{equation}
It follows from the fact that $C_n(x)$ is a polynomial of degree $n$ in
$x$ that $F(y;\alpha)$ is a 
polynomial of degree $n$ in $y$ (thus $y$ can
therefore be regarded as a continuous variable). We can therefore
write
\begin{equation}\label{a2.2}
F(y;\alpha) = \sum_{j=0}^n \kappa_{nj} C_j(y).
\end{equation}
To calculate the $\kappa_{nj}$ we first multiply both sides by
$q^{y/2} \alpha^{y/2} $ and subtract 
the same equation with
$y$ replaced by $y-1$ to obtain
\begin{equation}\label{a2.3}
C_n(y) = - F(y;\alpha) + q^{-1/2} \alpha^{-1/2} F(y-1;\alpha).
\end{equation}
To proceed further we note, as can be verified
from the definitions (\ref{3.1}) and (\ref{3.3}), that
\begin{eqnarray}
x C_n(x-1) & = & C_{n+1}(x) - {q \over q - 1} (2n+1) C_n(x) +
{n^2 q^2 \over (q-1)^2 } C_{n-1}(x), \nonumber \\
x C_n(x)  & = & C_{n+1}(x) - {n q + n + q \over q - 1} C_n(x) +
{n^2 q \over (q-1)^2 } C_{n-1}(x) \label{3C}
\end{eqnarray}
(the second of these is the three term recurrence, the general structure of
which holds for all sequences of orthogonal polynomials).
Thus with (\ref{a2.2}) substituted in (\ref{a2.3}), we can multiply both sides
of the equation by $y$ and use the equations to get an equation
involving only the linearly independent functions $C_{n+1}(y), C_n(y), \dots,
C_{0}(y)$ in the variable $y$. Equating coefficients of these functions
shows
\begin{equation}\label{a2.4}
\kappa_{nn} = \frac{\sqrt{q}}{\sqrt{\alpha} - \sqrt{q}}, \ \
\kappa_{nk} = \frac{n!}{k!} \frac{\sqrt{\alpha q}}{(\sqrt{\alpha}
- \sqrt{q})^2} \left(\frac{\sqrt{q} ( 1 - \sqrt{\alpha q} )}{
(1 - q)(\sqrt{\alpha} - \sqrt{q})} \right)^{n-k-1}, \ \ k < n.
\end{equation}

We will now make use of (\ref{a2.2}) and (\ref{a2.4}) to evaluate
(\ref{3.9}). Substituting (\ref{a2.2}) in (\ref{2.30a}), multiplying both
sides by $q^y C_m(y)$ and summing over $y$, making use of (\ref{ortho1}) on
the right hand side, we see that
\begin{equation}
\sum_{y=0}^{\infty} q^{y/2} \alpha^{-y/2} C_m(y) \sum_{x=y+1}^{\infty}
q^{x/2} \alpha^{x/2} C_n(x) = \left\{ \begin{array}{ll} h_m \kappa_{nm},
& m \leq n, \\ 0, & m > n. \end{array} \right.
\end{equation}
Recalling (\ref{antib}), and making use of the explicit formulas (\ref{3.6})
and (\ref{a2.4}), (\ref{3.15}) follows.
\hfill $\square$

\medskip
In general \cite[Prop.~6]{FR02a}, 
if the skew product has a factorization (\ref{3.15}) for monic
polynomials $C_n(x)$, $n=0,1,\dots$, the corresponding monic skew
orthogonal polynomials $R_k(x)$ can be written as a series in
$\{C_n(x)\}_{n=0,1,\dots,k}$ for explicit coefficients involving the
$a_j$'s and $b_j$'s, while $r_n = a_{2n} b_{2n+1}$.
These facts can be seen from the determinant formulas
(\ref{R1}) and (\ref{R2}). We thus have the following result.

\begin{prop}
Let $C_n(x)$ be specified by (\ref{3.3}). The monic skew orthogonal
polynomials with respect to (\ref{antib}) are given in terms of
these polynomials by
\begin{eqnarray}
R_{2n}(x) & = & C_{2 n}(x)
+ \sum_{k=0}^{n-1} \frac{(2 n)!}{(2 k)!}
\frac{q^{n-k}}{(1-q)^{2 n - 2 k}} C_{2 k}(x) \nonumber \\
&  & - \frac{\sqrt{\alpha} - \sqrt{q}}{1 - \sqrt{\alpha q}}
\sum_{k=0}^{n-1} \frac{(2 n)!}{(2 k+1)!}
\frac{q^{n-k-(1/2)}}{(1-q)^{2 n - 2 k-1}} C_{2 k+1}(x),
\label{reven}
\end{eqnarray}
\begin{equation}
R_{2 n + 1}(x) = C_{2 n + 1}(x) - \frac{\sqrt{\alpha}-\sqrt{q}}{1 -
\sqrt{\alpha q}} \frac{\sqrt{q}}{1-q} (2 n + 1) C_{2 n}(x),
\label{rodd}
\end{equation}
while the corresponding normalization has the explicit value
\begin{equation}\label{rn}
r_n =
\frac{(2 n)! (2 n + 1)! \sqrt{\alpha} q^{2 n + (1/2)}}{(1 - q)^{
4 n + 1} (1 - \sqrt{\alpha q})^2}.
\end{equation}
\end{prop}

\noindent
Proof. \quad These results follow from the general formulas of
\cite{FR02a}. The only point which
requires clarification is that we have chosen for the
arbitrary constant $c_n$ in (\ref{R2})
$$
c_n = - \frac{\sqrt{\alpha} - \sqrt{q}}{1 - \sqrt{\alpha q}}
\frac{\sqrt{q}}{1 - q} (2 n + 1),
$$
which gives
the simplest expression for $R_{2 n + 1}(x)$
in the basis $\{C_j(x)\}$.
\hfill $\square$

\medskip

We note that
following the procedure used in \cite{FNH98},
the triangular structures (\ref{reven}) and (\ref{rodd}) can be inverted to
give
\begin{eqnarray}\label{4.24}
& & C_{2 n}(x) = R_{2 n}(x) + \frac{\gamma^{2 n} (2 n)! q^n}{(1 - q)^{2 n}}
\nonumber \\
&& \quad \times 
\left[ \left( 1 - \frac{1}{\gamma^2} \right) \sum_{k=0}^{n-1}
\frac{(1-q)^{2 k}}{\gamma^{2 k} (2 k)! q^k} R_{2 k}(x)
 +  \sum_{k=0}^{n-1} \frac{(1-q)^{2 k+1}}{\gamma^{2 k+1} (2 k+1)! q^{k+(1/2)}}
R_{2 k+1}(x) \right] \nonumber \\
\end{eqnarray}
and
\begin{eqnarray}\label{4.25}
C_{2 n+1}(x)  & = & \gamma \frac{q^{1/2}}{1 - q} (2 n + 1) R_{2 n}(x)
 +  \frac{\gamma^{2 n+1} (2 n+1)! q^{n+(1/2)}}{(1 - q)^{2 n+1}}
\nonumber \\
&& \times
\left[ \left( 1 - \frac{1}{\gamma^2} \right) \sum_{k=0}^{n-1}
\frac{(1-q)^{2 k}}{\gamma^{2 k} (2 k)! q^k} R_{2 k}(x) +
\sum_{k=0}^n \frac{(1-q)^{2 k+1}}{\gamma^{2 k+1} (2 k+1)! q^{k+(1/2)}}
R_{2 k+1}(x) \right], \nonumber \\
\end{eqnarray}
where
\begin{equation}
\gamma = \frac{\sqrt{\alpha}-\sqrt{q}}{1 - \sqrt{\alpha q}}.
\end{equation} 

\subsection{A summation formula}
Let us write
\begin{equation}\label{bnn}
C_n(x) = \sum_{j=0}^n \beta_{nj}
R_j(x), \ \ \beta_{nn} = 1.
\end{equation}
We see from the explicit formulas (\ref{4.24}), (\ref{4.25}) that the
coefficients have the factorization property
\begin{equation}\label{5.24a}
\beta_{nj} = \alpha_n \gamma_j
\end{equation}
for certain $\alpha_n$, $\gamma_j$. According to \cite[Eq.~(3.42)]{NF}, under such
a circumstance the matrix element $S(x,y)$ can be summed.

\begin{prop}
The expression (\ref{sxy}) can be simplified to read
\begin{eqnarray}
 S(x,y) & = & q^{(x+y)/2} \frac{1}{h_{M-1}}
\frac{C_M(x) C_{M-1}(y) - C_{M-1}(x) C_M(y)}{x-y} 
+ q^{y/2} \frac{1}{r_{(M-2)/2}} \frac{1-q}{M \gamma q^{1/2}}
\nonumber \\ && \times
\left( \Phi_{M-2}(x) - q^{x/2} \frac{r_{(M-2)/2}}{h_{M-1}}
C_{M-1}(x) \right) \left( C_M(y) - R_M(y) \right).
\label{sxyr}
\end{eqnarray}
\end{prop}

\noindent
Proof. \quad We know from general formulas \cite{NF} that with
$\{\beta_{nj} \}_{j=0,\dots,n}$ the lower triangular transition matrix
for the change of variables from the monic skew orthogonal polynomials
$\{ R_j(x) \}_{j=0,1,\dots}$ to the monic orthogonal polynomials
$\{ C_j(x) \}_{j=0,1,\dots}$ as implied by (\ref{bnn}), the quantity
(\ref{phj}) permits the expansions
\begin{eqnarray}
\Phi_{2k-1}(x) & = & -
q^{x/2} \sum_{\nu = 2 k - 2}^{\infty} \frac{C_{\nu}(x)}{h_{\nu}}
\beta_{\nu \ 2 k - 2} \, r_{k-1}, \nonumber \\
\Phi_{2k-2}(x) & = & q^{x/2} \sum_{\nu = 2 k - 1}^{\infty} \frac{C_{\nu}(x)}
{h_{\nu}} \beta_{\nu \ 2 k - 1} \, r_{k-1}.
\label{phi}
\end{eqnarray}
Further, we know that these formulas substituted into (\ref{sxy}) give
\begin{equation}
S(x,y) = q^{(x+y)/2} \sum_{\nu = 0}^{M-1}
 \frac{C_{\nu}(x) C_{\nu}(y)}{h_{\nu}} +
q^{(x+y)/2} \sum_{\nu = M}^{\infty} \sum_{k=0}^{M-1}
\frac{C_{\nu}(x)}{h_{\nu}} \beta_{\nu k} R_k(y).
\label{sxys}
\end{equation}
According to the Christoffel-Darboux formula from the theory of orthogonal
polynomials, the first summation gives the first term on the 
right hand side of
(\ref{sxyr}). Moreover, making use of the factorization (\ref{5.24a}) shows
$$
\sum_{k=0}^{M-1} \beta_{\nu k} R_k(y) = \frac{\gamma^{\nu-M}
q^{(\nu-M)/2}}{(1 - q)^{\nu-M}} \frac{\nu!}{M!} \left( C_M(y) - R_M(y)
\right).
$$
After further use of (\ref{phi}) the second term of (\ref{sxyr})
results. \hfill $\square$

\section{Poissonized Random Involutions}
\setcounter{equation}{0}

\subsection{Poissonization}

Essential to the rationale underlying the introduction of the geometrical model
in \S \ref{spdf}
is its reduction
in an appropriate limit
to the Poissonization of random involutions. To see how this comes
about, the procedure for the Poissonization of random involutions
given in \cite{BR02a} must be revised.

Consider involutions of $\{1,2,\dots,N\}$, and specify the number of
fixed points therein to be $m$. Catalogue the involutions according to
$(\lambda_N^{(k)} =: \lambda_k)_{k=1,\dots,N}$ where $\lambda_N^{(k)}$ is
specified by (\ref{1.b}). The Robinson-Schensted-Knuth correspondence
(see e.g.~\cite{Fu97}) tells us that there is a bijection between all
involutions with a given value of $\lambda = (\lambda_1,\lambda_2,\dots,
\lambda_N)$ and standard tableaux of shape $\lambda$ and content $N$.
Furthermore, the number of fixed points $m$ is equal to
$\sum_{j=1}^N (-1)^{j-1} \lambda_j$, and so restricts the permissible
$\lambda$.

Let the number of standard tableaux of shape $\lambda$ and content
$N$ be denoted $f_\lambda$. With $\ell(\lambda)$ denoting the number of
non-zero parts of $\lambda$, and
\begin{equation}\label{2.xc}
V_p(\lambda) := \prod_{1 \leq j < l  \leq p}(\lambda_j - \lambda_l + l - j),
\qquad
W_p(\lambda) = \prod_{j=1}^p \frac{1}{(\lambda_j + p - j)!},
\end{equation}
it is known that \cite{Fu97}
\begin{equation}
f_{\lambda} = N! V_{l(\lambda)}(\lambda) W_{l(\lambda)}.
\end{equation}
Further, let $t_{n,m}$ denote the number of involutions of
$\{1,2,\dots,N\}$ with $m$ fixed points and $n$ 2-cycles. We have
\begin{equation}\label{tnm}
t_{n,m} = \sum_{\lambda \in  \Omega_{2n+m,m}^{(\infty)}}  f_{\lambda}
= {(2n+m)! \over 2^n n! m!}, 
\end{equation}
where $\Omega_{2n+m,m}^{(\infty)}$ is the set of
partitions with
\begin{equation}\label{3.2a}
\sum_{j=1}^{\ell(\lambda)}
 \lambda_j := |\lambda| = 2n+m, \qquad \sum_{j=1}^{\ell(\lambda)} (-1)^{j-1} \lambda_j
= m.
\end{equation}

In terms of these quantities the probability that such an involution
corresponds to a standard tableau of shape $\lambda$ is then given by
$$
\mathbb Q_{n,m}(\lambda) = {f_\lambda \over t_{n,m} }
\chi_{\lambda \in \Omega_{2n+m,m}^{(\infty)}}.
$$
And the Poissonized form of this, in which $n$ and $m$ are regarded
as random variables from distinct Poisson distributions, is
\begin{eqnarray}\label{Q}
\mathbb Q^{(z_1,z_2)}(\lambda)
& := & {e^{-z_1 - z_2} z_1^n z_2^m \over n! m!}
{1 \over t_{n,m} } f_\lambda 
\nonumber \\ 
& = &  e^{-z_1 - z_2} z_1^n z_2^m {2^n \over (2n+m)!}
f_\lambda,
\end{eqnarray}
where the second equality follows upon using (\ref{tnm}), 
and $n$, $m$ are related to $\lambda$ by (\ref{3.2a}).

For future applications it is convenient to write this in terms of
$(N,m)$ rather than $(n,m)$. Setting
\begin{equation}\label{aQ}
z_1 = Q/2, \qquad  z_2 = \sqrt{Q\alpha}
\end{equation}
and with $|\lambda| = N$ defining
\begin{equation}\label{aP}
\mathbb P^{(Q,\alpha)}(\lambda) = e^{- \sqrt{\alpha Q} - Q/2}
{\alpha^{m/2} Q^{N/2} \over N!}
f_\lambda, 
\end{equation}
we see that
\begin{equation}\label{QP}
\mathbb Q^{(Q/2,(\alpha Q)^{1/2})}(\lambda) =
\mathbb P^{(Q,\alpha)}(\lambda).
\end{equation}
Also for future application, we note from the fact that
$ \mathbb P^{(Q,\alpha)}$ is a probability distribution 
that one has the summation 
\begin{equation}\label{3.b}
\sum_{N=0}^\infty
{Q^{N/2} \over N!} t_N^{(\alpha)} = e^{\sqrt{\alpha Q} + Q/2},
\end{equation}
where
\begin{equation}\label{3.1t}
 t_N^{(\alpha)} := \sum_{\lambda: |\lambda| = N}  
t_{n,m} \alpha^{m/2}
\end{equation}
(the case $\alpha = 1$ is given in \cite[\S 5.1.4, eq.~(42)]{Kn98}).

Consider now a function $g(\lambda)$ of the parts of a partition $\lambda$.
A concrete example of future use is
\begin{equation}\label{4.1a}
g(\lambda) = \prod_{j=1}^\infty (1 + u(\lambda_j - j) ),
\end{equation}
where $u(x)$ is any function which vanishes for $x < 0$. For such
functions we define the  Poissonized average as
\begin{equation}\label{4.g}
\langle g(\lambda) \rangle^{(Q,\alpha)}
=  
\sum_{\lambda }
g(\lambda) 
\mathbb P^{(Q,\alpha)}(\lambda).
\end{equation}
We want to show that this average results as a limit of the model of
\S \ref{spdf}. 

\begin{prop}\label{p5}
Let $g(\lambda)$ be a function of the parts of a partition $\lambda$
which furthermore satisfies
the bound
\begin{equation}\label{5.13}
|g(\lambda)| \le c^{\sum_{j=1}^{\ell(\lambda)} \lambda_j}
\end{equation}
for some $c > 0$. With $P_M(\lambda;q)$ given by (\ref{1.9P}) we have
\begin{equation}\label{5.e}
\lim_{M \to \infty}
\sum_{\lambda: \ell(\lambda) \le M} g(\lambda)
P_M(\lambda;Q/M^2) = \langle g(\lambda) \rangle^{(Q,\alpha)}.
\end{equation}
\end{prop}

\noindent
Proof. \quad We essentially follow Johansson 
\cite{KJ2}, who proved the analogous
result for Poissonized permutations. Using the explicit formula
(\ref{1.9P}), we see 
\begin{eqnarray}
\sum_{\lambda: \ell(\lambda) \le M} 
g(\lambda) P_M(\lambda;q)
& = &  (1 - \sqrt{\alpha q})^M (1 - q)^{M(M-1)/2}
\sum_{\lambda: \ell(\lambda) \le M} g(\lambda)
 \prod_{j=1}^M q^{\lambda_j/2}
\prod_{j=1}^{M} \alpha^{\sum_{j=1}^M(-1)^{j-1} \lambda_j} \nonumber \\
&& \qquad \times 
\prod_{1 \le j < k \le M}
{\lambda_j - \lambda_l + l - j \over l - j}  \nonumber \\
& = & (1 - \sqrt{\alpha q})^M (1 - q)^{M(M-1)/2}
\sum_{\lambda: \ell(\lambda) \le M}
q^{N/2} \alpha^{m/2}
g(\lambda) V_M(\lambda) \prod_{j=1}^{M-1} {1 \over j!}. \nonumber \\
\label{5.8}
\end{eqnarray}
With $\Omega_{N,m}^{(M)}$ denoting the set of partitions $\lambda$ with the properties
(\ref{3.2a}) but constrained so that $\ell(\lambda) \le M$, 
recalling (\ref{2.xc}) we see
\begin{eqnarray}\label{5.12}
\sum_{\lambda \in \Omega_{N,m}^{(M)}} g(\lambda) V_M(\lambda)
\prod_{j=1}^{M-1} {1 \over j!}
& = & \sum_{\lambda \in \Omega_{N,m}^{(M)}}
g(\lambda) \frac{V_{l(\lambda)}(\lambda)
W_{l(\lambda)}(\lambda)}{W_M(\lambda)}
\prod_{j=1}^{M-1} \frac{1}{j!}
\nonumber \\
& =  &  \sum_{\lambda \in \Omega_{N,m}^{(\infty)}} 
g(\lambda) \frac{f_{\lambda}}{N!} 
\prod_{j=1}^{M} \frac{(\lambda_j+M-j)!}{(M-j)!}.
\end{eqnarray}

From the assumed bound (\ref{5.13}) and the further bound
$$
{(\lambda_j + M - j)! \over M^{\lambda_j} (M-j)!} \le 1
$$
(a consequence of Stirling's formula), we see from (\ref{5.12})
that with $q$ replaced by $Q/M^2$ in (\ref{5.8}), the sum over
$\lambda: \ell(\lambda) \le M$
for $M \to \infty$ in the latter is itself bounded
by
$$
\sum_{\lambda}
{\tilde{Q}^{N/2} \alpha^{m/2} \over N!} t_{N,m},
\qquad \tilde{Q} = Qc^2.
$$
But according to (\ref{3.1t}) and (\ref{3.b}) this is summed as the
r.h.s.~of (\ref{3.b}) with $Q$ replaced by $\tilde{Q}$. Thus we have
a uniform bound in $N$ for the sum in (\ref{5.8}) with $q$ replaced
by $Q/M^2$. The limit $M \to \infty$ can therefore be taken term-by-term.
Doing so gives (\ref{5.e}).
\hfill $\square$

\smallskip
We remark that (\ref{4.g}) can be written
\begin{eqnarray}
\langle g(\lambda) \rangle^{(Q,\alpha)} & = &
e^{-\sqrt{\alpha Q} - Q/2} \sum_{N=0}^\infty
{Q^{N/2} \over N!} t_N^{(\alpha)} \langle g(\lambda) \rangle^{(\alpha)}_N,
\label{aN1} \\
\langle g(\lambda) \rangle^{(\alpha)}_N & := &
{1 \over t_N^{(\alpha)} }
\sum_{\lambda: |\lambda| = N} g(\lambda) \alpha^{m/2}
f_\lambda.
\label{aN}
\end{eqnarray}
Some insight into the role of $\alpha$ can be obtained by computing
the mean number of fixed points in the ensemble (\ref{aN}). From (\ref{4.g}),
the
identity (\ref{QP}), the constraint involving $m$ in (\ref{3.2a}), and the explicit
form (\ref{Q})
it is immediate that
$$
\Big \langle \sum_{j=1}^{\ell(\lambda)} (-1)^{j-1} \lambda_j \Big \rangle^{(Q,\alpha)}
= \sqrt{\alpha Q}.
$$
This in (\ref{aN1}) and use of (\ref{3.b}) implies
$$
\Big \langle \sum_{j=1}^{\ell(\lambda)} (-1)^{j-1} \lambda_j \Big \rangle^{(\alpha)}_N
= \sqrt{\alpha} {t_{N-1}^{(\alpha)} \over t_N^{(\alpha)} } N.
$$
Recalling the explicit formula (\ref{tnm}), the sum (\ref{3.1t}) can be
estimated to give $t_{N}^{(\alpha)} / t_{N-1}^{(\alpha)} \sim \sqrt{N}$
as $N \to \infty$ ($\alpha > 0$). Hence
$$
\Big \langle \sum_{j=1}^{\ell(\lambda)} (-1)^{j-1} \lambda_j \Big \rangle^{(\alpha)}_N
\mathop{\sim}\limits_{N \to \infty} \sqrt{\alpha N},
$$
which tells us that in the ensemble (\ref{aN}) the mean number of
fixed points is necessarily proportional to $\sqrt{N}$, with
proportionality constant $\sqrt{\alpha}$.

\subsection{A Poissonized average and limiting correlations}
According to the equations (\ref{6a.c})--(\ref{6a.2}), to calculate
the distribution function Pr$(h_1 \le a_1, \dots, h_l \le a_l)$, it is
sufficient to compute the average
\begin{equation}\label{dis}
\Big \langle \prod_{l=1}^M\Big (1 - \sum_{r=1}^k \xi_r \chi_{I_r}^{(l)}
\Big )
\Big \rangle_{\rm Psym}.
\end{equation}
Let us suppose each end point $a_j$ in the $\{I_r\}$ is a positive
integer when measured from $M$. We indicate this by the replacement
$I_r \mapsto M + I_r$. The quantity being averaged in (\ref{dis}) is now
of the form (\ref{4.1a}), and   (\ref{5.e}) tells us that
\begin{equation}\label{3.2ab}
\lim_{M \to \infty} 
\Big \langle \prod_{l=1}^M \Big (1 - \sum_{r=1}^k \xi_r \chi_{I_r}^{(l)}
\Big )
\Big \rangle_{\rm Psym} \Big |_{q = Q/M^2}
=
\Big \langle \prod_{l=1}^\infty
\Big (1 - \sum_{r=1}^k \xi_r \chi_{\lambda_l - l \in
I_r} \Big ) \Big \rangle^{(\alpha,Q)}.
\end{equation}

We know from the proof of Proposition \ref{p5} that the convergence is 
uniform in $M$, and can therefore be taken term-by-term in
(\ref{6a.2}) to give
\begin{eqnarray}\label{3.2ac}
\lefteqn{
\Big \langle \prod_{l=1}^\infty
\Big (1 - \sum_{r=1}^k \xi_r \chi_{\lambda_l - l \in
I_r} \Big ) \Big \rangle^{(\alpha,Q)} } \nonumber \\
&& =
1 + \sum_{p=1}^\infty {(-1)^p \over p!} \sum_{h_1,\dots,h_p = 0}^\infty
\prod_{i=1}^p \Big ( \sum_{r=1}^k \xi_r \chi_{I_r}^{(i)} \Big )
\lim_{M \to \infty}
\rho_{(p)}(h_1+M,\dots,h_p+M) \Big |_{q = Q/M^2}. \nonumber \\ 
\end{eqnarray}
Moreover, the uniform convergence tells us that the limit $M \to \infty$
can be taken term-by-term in the summations of the quantities
specifying $\rho_{(k)}$ (recall (\ref{sxy})--(\ref{phj})).
To calculate these limits requires an appropriate asymptotic
formula relating to $C_n(x)$.

\begin{prop}
For $M \to \infty$ we have
\begin{equation}\label{climit}
C_n(x+M) \Big |_{q = Q/M^2}
\sim  \frac{n! q^{(n-M-x)/2}}{(1-q)^n}
J_{-n+M+x}(2 \sqrt{Q})
\end{equation}
uniformly in $n \in \mathbb Z_{\ge 0}$.
\end{prop}

\noindent
Proof. \quad
Following \cite{KJ1},
we use the generating
function of the Meixner polynomials \cite{KOEKOEK}
\begin{equation}
\left( 1 - \frac{t}{q} \right)^x (1 - t)^{-x-c}
= \sum_{n=0}^{\infty} \frac{(c)_n}{n!} M_n(x;c,q) t^n,
\end{equation}
to derive the integral representation
\begin{equation}
M_n(x;c,q) =  \frac{n!}{2 \pi i (c)_n}
\oint {\rm d}t \ t^{-n-1} \left( 1 - \frac{t}{q} \right)^x
(1 - t)^{-x-c},
\label{integral}
\end{equation}
where the path of integration encloses the origin anticlockwise.
It follows from  (\ref{integral}) that
\begin{eqnarray}
& & C_n(x+M)|_{q = Q/M^2} \nonumber \\ & = & \frac{n! q^n}{2 \pi i (q-1)^n}
\oint {\rm d}t \ t^{-n-1} \left( 1 - \frac{t}{q} \right)^{x+N}
(1 - t)^{-x-N-1} \nonumber \\
& = & \frac{n! q^{(n-M-x)/2}}{(1-q)^n} \frac{(-1)^{-n+M+x}}{2 \pi}
\int_{-\pi}^{\pi} {\rm d}\theta \ z^{-n+M+x} \left( 1 - \frac{\sqrt{Q}}{M z}
\right)^{x+M} \left(1 - \frac{\sqrt{Q} z}{M} \right)^{-x-M-1}
\label{cintegral}
\end{eqnarray}
with $z = r {\rm e}^{i \theta}$, $r>0$. For large $M$ the integral has the
leading form
$$
\int_{-\pi}^\pi {\rm d}\theta \,  \ z^{-n+M+x}  \exp\left\{ \sqrt{Q}\left(
z - \frac{1}{z} \right)\right\} = (-1)^{-n+M+x} 2 \pi
J_{-n+M+x}(2 \sqrt{Q})
$$
uniform in $n \in \mathbb Z$, where $J_k(x)$ denotes the Bessel function,
and the result follows. \hfill $\square$

\medskip
Asymptotic formulas for $\{R_j(x)\}$, $\{\Phi_j(x)\}$ appearing in the
matrix elements of (\ref{2.11}) can now be obtained.

\begin{cor}\label{c2}
Let $M$ be even. We have
\begin{eqnarray}
R_{M-2n}(x+M) \Big |_{q = Q/M^2}
& \sim &   (M - 2n)! q^{-(2n+x)/2}\Big (
\sum_{l=0}^{M/2-n}  J_{2 n +2l + x}(2 \sqrt{Q})
\nonumber \\
&&  
- \sqrt{\alpha}
\sum_{l=1}^{M/2-n} J_{2 n +2l -1 + x}(2 \sqrt{Q}) \Big ),
\label{rlimit}
\end{eqnarray} 
\begin{equation}\label{rlimit1}
R_{M-2n+1}(x+M) \Big |_{q = Q/M^2}
 \sim (M - 2n + 1)! q^{-(2n+x-1)/2} \Big (
J_{2n - 1+x}(2 \sqrt{Q}) - \sqrt{\alpha} J_{2n +x}(2 \sqrt{Q}) \Big ).
\end{equation}
\begin{equation}
\Phi_{M-2n}(x+M) \Big |_{q = Q/M^2}
 \sim  \alpha^n (M-2n)! q^{M-2n)/2}
\sum_{\nu = -2n+1}^\infty \alpha^{\nu/2} J_{-\nu + x}(2 \sqrt{Q}),
\label{philimit}
\end{equation}
\begin{eqnarray}
 \Phi_{M-2n+1}(x+M) \Big |_{q = Q/N^2} & \sim & \alpha^{n+1/2}
(M-2n+1)! q^{M-2n+1)/2}\nonumber \\
&&  \times
\sum_{\nu = -2n}^\infty \alpha^{\nu/2} \Big ( - J_{-\nu + x}(2 \sqrt{Q})
+ J_{\nu - 2 + x}(2 \sqrt{Q}) \Big ).
\label{philimit1}
\end{eqnarray}
\end{cor}

\noindent
Proof. \quad 
The expansions (\ref{rlimit}) and (\ref{rlimit1}) follow immediately upon
use of (\ref{climit}) in (\ref{reven}) and (\ref{rodd}). In relation to
$\{\Phi_j(x)\}$ we read off explicit formulas for the $\beta_{nj}$ as
defined in (\ref{bnn}) from (\ref{4.24}), (\ref{4.25}) and substitute in
(\ref{phi}) to obtain
\begin{eqnarray}
\Phi_{2 n}(x) & = & q^{x/2} \frac{\sqrt{\alpha}}{(1 -
\sqrt{\alpha q})^2} \frac{(2 n)!
q^n}{\gamma^{2 n + 1} (1 - q)^{2 n}}
\sum_{\nu=2 n + 1}^{\infty}
\frac{(1-q)^{\nu+1} \gamma^{\nu}}{\nu! q^{\nu/2}}
C_{\nu}(x), \nonumber \\
\Phi_{2 n + 1}(x) & = & q^{x/2} \frac{\sqrt{\alpha}}{(1 -
\sqrt{\alpha q})^2} \frac{(2 n+1)!
q^{n+(1/2)}}{\gamma^{2 n} (1 - q)^{2 n + 1}} \nonumber \\
& & \times
\sum_{\nu=2 n}^{\infty}
\frac{(1-q)^{\nu+1} \gamma^{\nu}}{\nu! q^{\nu/2}}
\left( -  C_{\nu}(x) + \frac{(1-q)^2}{(\nu+2)(\nu+1)q}
C_{\nu+2}(x) \right).
\label{phir}
\end{eqnarray}
Use of (\ref{climit}) now gives (\ref{philimit}), (\ref{philimit1}).
\hfill $\square$

\begin{prop}\label{p7}
We have
\begin{eqnarray}\label{3.33}
\rho_{(k)}^{(Q,\alpha)}(h_1,\dots,h_k) & := & 
\lim_{M \to \infty}
\rho_{(k)}(h_1+M,\dots,h_k+M) \Big |_{q = Q/M^2}
\nonumber \\
& = & 
{\rm qdet} \left [ \begin{array}{cc} \bar{S}(h_j,h_l) & \bar{I}(h_j,h_l) \\
\bar{D}(h_j,h_l) & \bar{S}(h_l,h_j) \end{array} \right ]_{j,l=1,\dots,k}, 
\end{eqnarray}
where
\begin{eqnarray}
\bar{S}(x,y) & =   &
 \frac{\sqrt{Q}}{x-y} \left(
J_x(2 \sqrt{Q}) J_{y+1}(2 \sqrt{Q}) -
J_y(2 \sqrt{Q}) J_{x+1}(2 \sqrt{Q}) \right) \nonumber \\
&  & \qquad - \sum_{j=0}^{\infty} \alpha^{j/2} J_{-j+x}(2 \sqrt{Q})
\sum_{l=0}^{\infty} \left( J_{2 l + 2 + y}(2 \sqrt{Q}) -
\sqrt{\alpha} J_{2 l + 1 + y}(2 \sqrt{Q}) \right), 
\label{slimit}
\end{eqnarray}
\begin{eqnarray}\label{ilimit}
 \bar{I}(x,y) & = &    - \sum_{n=1}^\infty \sum_{k=1}^\infty  \alpha^{k/2}
\nonumber \\
&  &  
\Big ( J_{n+x}(2 \sqrt{Q}) J_{n-k+y}(2 \sqrt{Q}) -
J_{n+y}(2 \sqrt{Q}) J_{n-k+x}(2 \sqrt{Q}) \Big ) + \epsilon (x,y),
\end{eqnarray}
\begin{eqnarray}\label{dlimit}
& & {\bar D}(x,y) \nonumber \\
&  &  \quad =
{1 \over \sqrt{\alpha}} \sum_{l=1}^\infty \sum_{j=1}^l
\Big ( J_{2l+x}(2 \sqrt{Q}) - \sqrt{\alpha} J_{2l+x+1}(2 \sqrt{Q}) \Big )
\Big ( J_{2j+y-1}(2 \sqrt{Q}) - \sqrt{\alpha} J_{2j+y}(2 \sqrt{Q}) \Big )
 \nonumber \\
&& \qquad
- (x \leftrightarrow y). 
\end{eqnarray}
\end{prop}

\noindent
Proof. \quad Using the results of Corollary \ref{c2}, from the form
of the matrix elements (\ref{sxy})--(\ref{phj}) we compute
$$
\bar{S}(x,y) = \lim_{M \to \infty}S(x+M,y+M) \Big |_{q = Q/M^2} 
$$
and similarly for $\bar{I}(x,y)$, $\bar{D}(x,y)$.
\hfill $\square$

\medskip
We note that Bessel function identities can be used to verify that the
first term in (\ref{slimit}) can be written in a denominator free from
\cite[Prop.~2.9]{BOO}. This allows (\ref{slimit}) to be replaced by 
\begin{eqnarray}
\bar{S}(x,y) & = & \sum_{n=1}^\infty J_{n+x}(2 \sqrt{Q})
J_{n+y}(2 \sqrt{Q})
\nonumber \\
&  & \qquad - \sum_{j=0}^{\infty} \alpha^{j/2} J_{-j+x}(2 \sqrt{Q})
\sum_{l=0}^{\infty} \left( J_{2 l + 2 + y}(2 \sqrt{Q}) -
\sqrt{\alpha} J_{2 l + 1 + y}(2 \sqrt{Q}) \right).
\label{slimit1}
\end{eqnarray}
\hfill $\square$

\medskip

The correlations (\ref{3.33}) relate to the measure on diagrams of partitions
(\ref{aP}).
A closely related measure, introduced in
\cite{BR02a}, is
\begin{equation}\label{bP}
\tilde{\mathbb P}^{(Q,\beta)}(\lambda) = e^{- \sqrt{\beta Q} - Q/2}
{\beta^{m/2} Q^{N/2} \over N!}
f_\lambda, 
\end{equation}
where $N = |\lambda|$ and $m = \sum_{j=1}^{\ell(\lambda)} (-1)^{j-1} \lambda_j'$ with
$\lambda_j'$ denoting the length of the $j$th column of the diagram of $\lambda$. 
This relates to decreasing subsequences in the involution. In the case
$\beta = 0$ only rows of even length are permitted. For this model the
corresponding $k$-point correlations have recently been computed by
Ferrari \cite[Lemma 5.2]{Fe04}, in the context of the polynuclear growth model 
from a flat substrate. The results obtained have a very similar
structure to that exhibited in Proposition \ref{p7}. For example, the density
(one-point correlation) is computed for (\ref{bP}) with $\beta = 0$ as
\begin{equation}
\tilde{\rho}_{(1)}^{(Q,0)}(x) =
\sum_{n=1}^\infty \Big ( J_{n+x}(2 \sqrt{Q})\Big )^2 -
J_{x+1}(2  \sqrt{Q}) \Big ( \sum_{m=1}^\infty J_{x+2m-1}(2 \sqrt{Q}) -
{(1 - (-1)^x) \over 2} \Big ),
\end{equation}
while for  (\ref{aP}) with $\alpha = 0$ Proposition \ref{p7} and (\ref{slimit1}) give
\begin{equation}
{\rho}_{(1)}^{(Q,0)}(x) =
\sum_{n=1}^\infty \Big ( J_{n+x}(2 \sqrt{Q})\Big )^2 -
J_{x}(2  \sqrt{Q})   \sum_{m=1}^\infty J_{x+2m}(2 \sqrt{Q}). 
\end{equation}
We also make mention that in the works \cite{Ba04,FR02a,IS04} the
parameter dependent correlations have recently been calculated for
other models generalizing the geometric picture of random involutions.

\section{Asymptotic Correlation Functions}
\setcounter{equation}{0}

\subsection{De-Poissonization}
Following the pioneering work \cite{Jo98}, the  following
de-Poissonization lemma was derived in \cite{BR02b}.

\begin{lemma}\label{lem2}
Let $d > 0$ be an arbitrary positive real number, and set
\begin{eqnarray*}
\mu_n^{(d)} & = & n + (2 \sqrt{d+1} + 1) \sqrt{n \log n}, \nonumber \\
\nu_n^{(d)} & = & n - (2 \sqrt{d+1} + 1) \sqrt{n \log n}.
\end{eqnarray*}
Let $q = \{ q_{n_1,n_2} \}_{n_1,n_2 \ge 0}$ be such that
$$
q_{n_1,n_2} \ge q_{n_1+1,n_2}, \quad
q_{n_1,n_2} \ge q_{n_1,n_2+1}, \quad
0 \le q_{n_1,n_2} \le 1
$$
and define
$$
\phi(\lambda_1,\lambda_2) = e^{- \lambda_1 - \lambda_2}
\sum_{n_1,n_2 \ge 0} q_{n_1,n_2} {\lambda_1^{n_1} \lambda_2^{n_2}  \over
n_1! n_2!}.
$$
Then there exists constants $C$ and $n_0$ such that for all $n_1,n_2 \ge
n_0$
\begin{equation}
\phi(\mu_{n_1}^{(d)}, \mu_{n_2}^{(d)}) - C(n_1^{-d} + n_2^{-d})
\le q_{n_1,n_2} \le \phi(\nu_{n_1}^{(d)}, \nu_{n_2}^{(d)}) + 
C(n_1^{-d} + n_2^{-d}).
\end{equation}
\end{lemma}

In using this lemma, we take $q_{n,m}$ as the joint distribution of
$\{\lambda_{n,m}^{(j)} \}_{j=1,\dots,l}$ for random involutions of
$\{1,\dots,N\}$, $N = 2n+m$ with $n$ 2-cycles and $m$ fixed points,
and thus
$$
q_{n,m} = {\rm Pr}(\lambda_{n,m}^{(1)} \le a_1, \dots ,
\lambda_{n,m}^{(l)} \le a_l).
$$
We know from the argument of \cite[Lemma 7.5]{BR02b} that this
satisfies the monotonicity conditions required for the validity of
Lemma \ref{lem2}. We further choose
$$
a_j = 2 \sqrt{Q} + Q^{1/6} s_j, \qquad
n_1 = [Q/2], \qquad
n_2 = [\sqrt{Q} - 2 w Q^{1/3}].
$$
The lemma then gives
\begin{eqnarray}\label{4.2}
&& \lim_{Q \to \infty} e^{-z_1 - z_2} \sum_{n,m \ge 0}
{z_1^n z_2^m \over n! m!}
{\rm Pr}(\lambda_{n,m}^{(1)} \le a_1, \dots ,\lambda_{n,m}^{(l)} \le a_l)
\Big |_{\begin{array}{l} \scriptstyle a_j = 2 \sqrt{Q} + Q^{1/6} s_j \\
\scriptstyle z_1 = Q/2 \\ 
\scriptstyle z_2 = \sqrt{Q} - 2 w Q^{1/3} \end{array} }
\nonumber \\
&& \quad =
\lim_{N \to \infty}
{\rm Pr}(\lambda_{n,m}^{(1)} \le a_1, \dots ,\lambda_{n,m}^{(l)}
\le a_l)
\Big |_{\begin{array}{l} \scriptstyle a_j = 2 \sqrt{Q} + Q^{1/6} s_j \\
\scriptstyle
n = [Q/2] \\ \scriptstyle m = [\sqrt{Q} - 2 w Q^{1/3}] \end{array} },
\end{eqnarray}
assuming the limit on the left hand side exists. Regarding the left hand
side, we know from (\ref{QP}) that it is equal to 
\begin{equation}\label{ch}
\lim_{Q \to \infty} {\rm Pr}^{(Q,\alpha)}(\lambda_1 \le a_1, \dots,
\lambda_l \le a_l) \Big |_{a_j = 2 \sqrt{Q} + Q^{1/6} s_j \atop
\sqrt{\alpha} = 1 - 2 w/Q^{1/6}},
\end{equation}
where Pr${}^{(Q,\alpha)}$ refers to the ensemble defined by the measure
(\ref{aP}).

As demonstrated in (\ref{6a.c})--(\ref{6a.2}), the probability in (\ref{ch}) is
fully determined by the right hand side of (\ref{3.2ab}). Using (\ref{3.2ab})
and (\ref{3.2ac}) this in turn can be written
\begin{eqnarray}\label{3.2ad}
\lefteqn{
\Big \langle \prod_{l=1}^\infty
\Big (1 - \sum_{r=1}^k \xi_r \chi_{\lambda_l - l \in
I_r} \Big ) \Big \rangle^{(\alpha,Q)} } \nonumber \\
&& =
1 + \sum_{p=1}^\infty {(-1)^p \over p!} \sum_{h_1,\dots,h_p = 0}^\infty
\prod_{i=1}^p \Big ( \sum_{r=1}^k \xi_r \chi_{I_r}^{(i)} \Big )
\rho_{(p)}^{(Q,\alpha)}(h_1,\dots,h_p). 
\end{eqnarray}
Therefore our task is to compute the $Q \to \infty$ limit of this
quantity with
\begin{equation}\label{mix}
I_r = (a_r,a_{r-1}), \quad a_0=\infty,  a_j = 2 \sqrt{Q} + Q^{1/6}
s_j \: \: (j=1,\dots,k), \quad \sqrt{\alpha} = 1 - 2w/Q^{1/6}.
\end{equation}

\subsection{The limit $Q \rightarrow \infty$}
We will make use of a lemma due to Soshnikov \cite[Lemma 2]{So02},
in a form close to that used in \cite{BF03}.

\begin{prop}\label{So}
Consider a sequence of point processes labelled by a parameter $L$.
Suppose that for $L \to \infty$, and after the linear scaling
$x_j \mapsto A_L(x_j) = \alpha_L + a_L x_j$ 
of each of the coordinates, the sequence
approaches a limit point process with correlations $\{\rho_k\}_{k=1,2,\dots}$
such that
\begin{equation}\label{zp}
\int_{-\infty}^\infty dx_1 \dots \int_{-\infty}^\infty dx_k \,
\prod_{i=1}^k \Big ( \sum_{r=1}^l \xi_r \chi_{I_r}^{(i)} \Big ) 
\rho_{(k)}(x_1,\dots,x_k)
= o(k!),
\end{equation}
and suppose furthermore that for each $r=1,\dots,l$
\begin{equation}\label{zp1}
\lim_{L \to \infty}  a_L^k \int_{A_L(I_r)} dx_1 \cdots
 \int_{A_L(I_r)} dx_k \, \rho_{(k)}^{(L)}(A_L(x_1),\dots,A_L(x_k))
= \int_{I_r} dx_1 \cdots  \int_{I_r} dx_k \, \rho_{(k)}(x_1,\dots,x_k).
\end{equation}
One then has
\begin{eqnarray}\label{4.8}
&& \lim_{L \to \infty}
\Big \langle \prod_i \Big ( 1 - \sum_{r=1}^l \xi_r \chi_{A_L(I_r)}^{(i)}
\Big ) \Big \rangle^{(L)} \nonumber \\ &&
=
1 + \sum_{k=1}^\infty {(-1)^k \over k!}
\int_{-\infty}^\infty dx_1 \cdots \int_{-\infty}^\infty dx_k \,
\prod_{i=1}^k \Big ( \sum_{r=1}^l \xi_r \chi_{I_r}^{(i)} \Big )
\rho_{(k)}(x_1,\dots,x_k).
\end{eqnarray}
\end{prop}

\noindent
The essential idea behind this result is that the condition
(\ref{zp}) implies the moment problem for the number of particles
in $\{I_r\}$ is definite, and thus the convergence of moments,
which are integrals of correlation functions, implies the convergence
of distributions.

In using Proposition \ref{So}, the sequence of point processes will be
those specified by the correlations (\ref{3.33}) and thus labelled
by the parameter $Q$. The limiting point
process is that specified by the correlations appearing on the
r.h.s.~of (\ref{zp1}). These correlations we will calculate to be
(\ref{1.11}). 
To deduce (\ref{4.8}) we must first
verify (\ref{zp}) for the correlations (\ref{1.11}).

\begin{lemma}
There exists a positive number $M$ such that
\begin{equation}\label{4.8a}
\rho_{(k)}^{\rm scaled}(x_1,\dots,x_k) \le e^{-(x_1+\cdots+x_k)}
k^{k/2} M^k
\end{equation}
for all $x_i \in [s_0,\infty)$ $(i=1,\dots,k)$, were $s_0$ is arbitrary
but fixed.
\end{lemma}

\noindent
Proof. \quad For $X \to \infty$ we know that
${\rm Ai}(X) = O(e^{-2X^{3/2}/3})$. Recalling (\ref{1.2so}), we see
from the explicit formulas (\ref{1.12}) that
\begin{eqnarray}\label{st}
&&f^{22}(X,Y) = f^{11}(Y,X) = O(e^{-2Y^{3/2}/3}) O(e^{uX/2}), \quad
f^{21}(X,Y) = O(e^{uX/2})O(e^{uY/2}), \nonumber \\
&&
f^{12}(X,Y) = O(e^{-2X^{3/2}/3})O(e^{-2Y^{3/2}/3}).
\end{eqnarray}
The qdet can be expanded \cite{Dy70}
\begin{equation}\label{4.9a}
{\rm qdet} \, [f(X_i,X_j)]_{i,j=1,\dots,k} =
\sum_{{\rm disjoint \: cycles} \atop \in S_k}
(- 1)^{k-l} \prod_1^l \Big ( f(X_a,X_b) f(X_b,X_c) \cdots
f(X_d,X_a) \Big )^{(0)},
\end{equation}
where the superscript $(0)$ denotes the operation $(1/2)$Tr, and $l$
denotes the number of disjoint cycles. We see from this that each term
consists of factors of the form $f^{11}(X,X)$, $f^{22}(X,X)$,
$f^{11}(X,Y) f^{22}(Y,X')$, or $f^{12}(X,Y) f^{21}(Y,X')$. The bounds
(\ref{st}) tell us that each term  is bounded by
$O(e^{-(x_1 + \cdots + x_k)})$ (of course a sharper bound can be given,
but this is sufficient for our purpose). 

We must also bound the
dependence on $k$. For this we note that in general the replacement
of the matrix elements of (\ref{1.11}) by
\begin{eqnarray}
&&
f^{11}(X,Y) \mapsto {a(X) \over a(Y)} f^{11}(X,Y), \quad
f^{12}(X,Y) \mapsto a(X) a(Y) f^{12}(X,Y), \nonumber \\
&& 
f^{21}(X,Y) \mapsto {1 \over a(X) a(Y)} f^{21}(X,Y)
\end{eqnarray}
leaves qdet unchanged. Choosing $a(X) = e^{X^{1+\mu}}$
$(0 < \mu \ll 1)$ we see from (\ref{st}) that each term is bounded for
$X,Y \to \infty$. According to Hadamard's bound for determinants, the
dependence on $k$ is therefore bounded by $k^{k/2} M^k$ for some
$M > 0$. \hfill $\square$

Because each $I_r$ is bounded from
below, the bound (\ref{4.8a}) establishes (\ref{zp}).
The remaining task is to verify (\ref{zp1}).
First we make note of an
alternative form of $f^{22}(X,Y)$.

\begin{lemma}
The formula for $f^{22}(X,Y)$ in (\ref{1.12}) can be written
\begin{eqnarray}\label{su.1}
f^{22}(X,Y) & = & K^{\rm soft}(X,Y) + {1 \over 2} \Big (
{\rm Ai}(Y) + {u \over 2} \int_Y^\infty {\rm Ai}(s) \, ds \Big )
\nonumber \\
&& \times e^{uX/2} \Big ( e^{-u^3/24} - \int_X^\infty
e^{-ut/2} {\rm Ai}(t) \, dt \Big ).
\end{eqnarray}
\end{lemma}

\noindent
Proof. \quad This follows by making use of the integral formula in
(\ref{1.2so}) for $K^{\rm soft}$, integration by parts, and use of the
Fourier transform
\begin{equation}\label{4.10a}
\int_{-\infty}^\infty e^{-yx} {\rm Ai}(x) \, dx = e^{-y^3/3}.
\end{equation}
\hfill $\square$

\medskip
It is precisely the form (\ref{su.1}) which appears in the asymptotic form
of $\bar{S}(x,y)$ relevant to verifying  (\ref{zp1}).

\begin{prop}\label{p12}
We have
\begin{equation}
Q^{1/6} \bar{S}(2\sqrt{Q} + Q^{1/6} X, 2\sqrt{Q} + Q^{1/6} Y)
\Big |_{\sqrt{\alpha} = 1 - 2w/Q^{1/6}} =
f^{22}(X,Y) \Big |_{u=-4w} + O(Q^{-1/6}) O(e^{-Y}).
\end{equation}
\end{prop}

\noindent
Proof. \quad Our main tool is the $\nu \to \infty$ asymptotic expansion
\begin{equation}\label{u2}
J_\nu(\nu - x(\nu/2)^{1/3}) \: \sim \:
\Big ( {2 \over \nu} \Big )^{1/3} {\rm Ai}(x) +
O({1 \over \nu}) O(e^{-x}),
\end{equation}
uniform in $x > x_0$ ($x_0$ fixed). As noted in \cite[Eq.~(4.11)]{BF03}, this is a
consequence of results due to Olver \cite{Ol74}. Using this with
$$
\nu = 2 \sqrt{Q} + Q^{1/6} X + n, \qquad
x = X + n Q^{-1/6}
$$
shows
\begin{eqnarray}\label{ss1}
&&
\sum_{n=1}^\infty J_{n+x}(2 \sqrt{Q})  J_{n+y}(2 \sqrt{Q}) 
\Big |_{x = 2 \sqrt{Q} + Q^{1/6} X \atop
y = 2 \sqrt{Q} + Q^{1/6} Y} \nonumber \\
&& \quad =
{1 \over Q^{1/3} }
\sum_{n=1}^\infty \Big ( {\rm Ai}(X+nQ^{-1/6}) +
O(Q^{-1/6}) O(e^{-(X+nQ^{-1/6})}) \Big ) \nonumber \\
&& \qquad \times
\Big ( {\rm Ai}(Y+nQ^{-1/6}) +
O(Q^{-1/6}) O(e^{-(Y+nQ^{-1/6})}) \Big ) \nonumber \\
&& \quad =
{1 \over Q^{1/6} } \Big (
\int_0^\infty {\rm Ai}(X+t) {\rm Ai}(Y+t) \, dt
+ O(Q^{-1/6}) O(e^{-(X+Y)}) \Big ),
\end{eqnarray}
where to obtain the final equality use has been made of the fact that
a Riemann sum appears in the previous equality.

The sum over $j$ in the second term of (\ref{slimit1}) is not suited to
the use of (\ref{u2}). To overcome this, we apply the generating function
expansion
$$
\sum_{n=-\infty}^\infty t^n J_n(z) = \exp \Big ( {z \over 2} (t
- {1 \over t} ) \Big )
$$
to conclude
$$
\sum_{j=0}^\infty \alpha^{j/2} J_{-j+x}(2 \sqrt{Q}) =
\alpha^{x/2} \Big ( e^{-\sqrt{Q}(\sqrt{\alpha} - 1/\sqrt{\alpha})}
- \sum_{j=1}^\infty \alpha^{-(j+x)/2} J_{j+x}(2\sqrt{Q}) \Big ).
$$
We can now use (\ref{u2}) to deduce the asymptotic formula 
\begin{eqnarray}\label{ss2}
&& \sum_{j=0}^\infty \alpha^{j/2} J_{-j+x}(2 \sqrt{Q})
\Big |_{x = 2 \sqrt{Q} + Q^{1/6} X \atop
\sqrt{\alpha} = 1 - 2w/Q^{1/6} } \nonumber \\
&& \quad =
e^{8 w^3/3} e^{-2wX} O(e^{-2Xw^2/Q^{1/6}})
- \int_0^\infty e^{2wt} {\rm Ai}(X+t) \, dt +
O(Q^{-1/6}) O(e^{-X}).
\end{eqnarray}

It remains to consider the sum over $l$ in the final line of
(\ref{slimit1}). For this we require in addition to (\ref{u2}) the
uniform asymptotic expansion \cite[Below (4.12)]{BF03} (again a consequence of
results in \cite{Ol74})
\begin{equation}\label{4.14a}
J_{\nu \pm 1}(\nu - x(\nu/2)^{1/3}) \: \sim \:
\Big ( {2 \over \nu} \Big )^{1/3} {\rm Ai}(x)
\pm \Big ( {2 \over \nu} \Big )^{2/3} {\rm Ai}'(x) + O( {1 \over \nu} )
O(e^{-x}).
\end{equation}
We find
\begin{eqnarray}\label{ss3}
&&\sum_{l=0}^\infty \Big ( J_{2l+2+y}(2 \sqrt{Q}) - \sqrt{\alpha}
J_{2l+1+y}(2 \sqrt{Q}) \Big ) \Big |_{\begin{array}{l}
\scriptstyle x = 2 \sqrt{Q} + Q^{1/6} X \\
\scriptstyle y = 2 \sqrt{Q} + Q^{1/6} Y \\
\scriptstyle \sqrt{\alpha} = 1 - 2w/Q^{1/6} \end{array} } \nonumber
\\&& \quad =
{1 \over 2 Q^{1/6} } \Big (
2 w \int_0^\infty {\rm Ai}(Y+t) \, dt - {\rm Ai}(Y)
+ O(Q^{-1/3}) O(e^{-Y}) \Big ).
\end{eqnarray}

Substituting (\ref{ss1}), (\ref{ss2}) and (\ref{ss3}) in
(\ref{slimit1}) gives the stated result.
\hfill $\square$

For the remaining matrix elements in (\ref{1.12}) the following forms
appear in the asymptotic form of $\bar{D}$ and $\bar{I}$.

\begin{lemma}
The formulas for $f^{12}$ and $f^{21}$ in (\ref{1.12}) can be written
\begin{eqnarray}
f^{12}(X,Y) & = & {1 \over 4} \Big ( {u \over 2} +
{\partial \over \partial X} \Big )
\Big ( {u \over 2} +
{\partial \over \partial Y} \Big ) \nonumber \\
&& \times
\int_0^\infty ds \int_s^\infty dt \, \Big (
{\rm Ai}(Y+s) {\rm Ai}(X+t) - {\rm Ai}(X+s) {\rm Ai}(Y+t) \Big ),
 \nonumber \\
f^{21}(X,Y) & = & - e^{u|X-Y|/2} {\rm sgn}(X-Y) \nonumber \\
&& + \Big ( \int_Y^\infty e^{u(Y-t)/2} K^{\rm soft}(X,t) \, dt
-  \int_X^\infty e^{u(X-t)/2} K^{\rm soft}(Y,t) \, dt \Big )
\nonumber \\
 & & - e^{u(Y-X)/2} e^{-u^3/24} \int_X^\infty e^{us/2}
{\rm Ai}(s) \, ds +
e^{u(X-Y)/2} e^{-u^3/24} \int_Y^\infty e^{us/2}
{\rm Ai}(s) \, ds. \nonumber \\
\end{eqnarray}
\end{lemma}

\noindent
Proof. \quad These follow straightforwardly upon using the integral
formula in (\ref{1.2so}) and the Fourier transform (\ref{4.10a}).
\hfill $\square$

\medskip
\begin{prop}\label{p13}
Let $0 < \mu < 1/2$ be fixed. We have
$$
Q^{1/3} \bar{D}(2 \sqrt{Q} + Q^{1/6} X, 2 \sqrt{Q} + Q^{1/6} Y)
\Big |_{\sqrt{\alpha} = 1 - 2w/Q^{1/6}}
= f^{12}(X,Y) \Big |_{u = -4w} +
O\Big ( {1 \over Q^{1/6}} \Big ) O(e^{-X^{1+\mu} - Y^{1+\mu}}),
$$
\begin{eqnarray*}
\lefteqn{\bar{I}(2 \sqrt{Q} + Q^{1/6} X, 2 \sqrt{Q} + Q^{1/6} Y) }
\nonumber \\
&& = f^{21}(X,Y) \Big  |_{u=-4w} +
e^{-2w |Y-X|} O( e^{2w^2 |Y-X|/Q^{1/6}} - 1) +
O\Big ( {1 \over Q^{1/6}} \Big ) \Big (  O(e^{-X}) +  O(e^{-Y}) \Big ).
\end{eqnarray*}
\end{prop}

\noindent
Proof. \quad The procedure used in the proof of Proposition
\ref{p12} suffices, with the error terms in (\ref{u2}) and
(\ref{4.14a}) sharpened from $O(e^{-X})$ to $O(e^{-X^{1+\mu}})$.
\hfill $\square$

\medskip
We remarked below the definition (\ref{4.9a}) of a quaternion
determinant that only specific combinations of the elements in
the underlying $2 \times 2$ matrix can occur. We can use this
fact together with the asymptotic expansions from
Propositions \ref{p12} and \ref{p13} to deduce an asymptotic
formula relating $\rho_{(k)}^{(Q,\alpha)}$ to
$\rho_{(k)}^{\rm scaled}$. This in turn implies
the $Q \to \infty$ limit of (\ref{3.2ad}).

\begin{cor}
For $Q \to \infty$ we have
\begin{eqnarray}\label{kt.1}
&& Q^{k/6} \rho_{(k)}^{(Q,\alpha)}(2 \sqrt{Q} + Q^{1/6} X_1,
\dots, 2 \sqrt{Q} + Q^{1/6} X_k) \Big |_{\sqrt{\alpha} = 1 - 2w/Q^{1/6}}
\nonumber \\
&& \qquad =
\rho_{(k)}^{\rm scaled}(X_1,\dots,X_k;-4w)
+ O\Big ( {1 \over Q^{1/6}} \Big )
O (e^{-(X_1 + \cdots + X_k)}). 
\end{eqnarray}
Consequently, with $I_r$ specified as in (\ref{mix}) and $s_0 := \infty$,
\begin{eqnarray}\label{kt.2}
&&\lim_{Q \to \infty}
\Big \langle \prod_{l=1}^\infty \Big ( 1 -
\sum_{r=1}^k \xi_r \chi_{\lambda_l - l \in I_r}
\Big ) \Big \rangle^{(\alpha,Q)}
\Big |_{\sqrt{\alpha} = 1 - 2w/Q^{1/6}} = 1  \nonumber \\
&&  + \sum_{p=1}^\infty {(-1)^p \over p!}
\int_{-\infty}^\infty dX_1 \dots \int_{-\infty}^\infty dX_p \,
\prod_{i=1}^p \Big (
\sum_{r=1}^k \xi_r \chi_{(s_r,s_{r-1})}^{(i)} \Big )
\rho_{(p)}^{\rm scaled}(X_1,\dots,X_p;-4w). \nonumber \\ 
\end{eqnarray}
\end{cor}

\noindent
Proof. \quad Consider first (\ref{kt.1}). The leading term is immediate.
Regarding the error bound, as implied by the
remark below (\ref{4.9a}), the elements in the $2 \times 2$ matrix of
(\ref{3.33}) appear as factors in (\ref{4.9a}) only in the 
combinations $\bar{S}(x,x)$, $\bar{S}(x,y) \bar{S}(y,x')$ or
$\bar{I}(x,y) \bar{D}(y,x')$. This, together with the individual
error bounds in Propositions \ref{p12} and \ref{p13}, imply
the error bound in (\ref{kt.1}).

We know from Proposition \ref{So} that to establish (\ref{kt.2}) it
suffices to establish (\ref{zp}) and (\ref{zp1}). As previously
remarked, the bound (\ref{4.8a}) implies (\ref{zp}). Thus it remains to
check that with $a_j$ as in (\ref{mix}),
\begin{eqnarray}\label{kt.3}
\lefteqn{
\lim_{Q \to \infty} Q^{k/6}
\sum_{h_1,\dots,h_p = 0}^\infty
\prod_{i=1}^p \chi^{(i)}_{(a_r,a_{r-1})}
\rho_{(p)}^{(Q,\alpha)}(h_1,\dots,h_p) \Big |_{\sqrt{\alpha} = 1 - 2w/
Q^{1/6}} } \nonumber \\
&& \qquad =
\int_{-\infty}^\infty dx_1 \cdots \int_{-\infty}^\infty dx_p \,
\prod_{i=1}^p \chi^{(i)}_{s_r,s_{r-1}}
\rho_{(p)}^{\rm scaled}(x_1,\dots,x_p;-4w).
\end{eqnarray}
For this we note
\begin{eqnarray*}
\lefteqn{\sum_{h_1,\dots,h_p = 0}^\infty
\prod_{i=1}^p \chi^{(i)}_{(2 \sqrt{Q} + Q^{1/6} s_r,
2 \sqrt{Q} + Q^{1/6} s_{r-1})}
\rho_{(p)}^{(Q,\alpha)}(h_1,\dots,h_p)} \nonumber \\
&& \quad
= \sum_{x_1,\dots,x_p: \atop
2\sqrt{Q} + Q^{1/6} x_j \in \mathbb Z}
\prod_{i=1}^p \chi^{(i)}_{(s_r,s_{r-1})}
\rho_{(p)}^{(Q,\alpha)}(2 \sqrt{Q} + Q^{1/6} x_1,\dots,
2 \sqrt{Q} + Q^{1/6} x_p).
\end{eqnarray*}
Recognising this as a Riemann sum, and substituting (\ref{kt.1})
gives (\ref{kt.3}).
\hfill $\square$.

\subsection{Proof of theorem 1}
According to (\ref{4.2}),  (\ref{ch}), (\ref{6a.c}) and (\ref{6a.1})
\begin{eqnarray}\label{f1f}
&&
\lim_{N \to \infty}
{\rm Pr}(\lambda_{n,m}^{(1)} \le a_1, \dots ,\lambda_{n,m}^{(l)}
\le a_l)
\Big |_{\begin{array}{l} \scriptstyle a_j = 2 \sqrt{Q} + Q^{1/6} s_j \\
\scriptstyle
n = [Q/2] \\ \scriptstyle m = [\sqrt{Q} - 2 w Q^{1/3}] \end{array} }
\nonumber \\&&
\quad = \lim_{Q \to \infty}
 {\rm Pr}^{(Q,\alpha)}(\lambda_1 \le a_1, \dots,
\lambda_l \le a_l) \Big |_{a_j = 2 \sqrt{Q} + Q^{1/6} s_j \atop
\sqrt{\alpha} = 1 - 2 w/Q^{1/6}} \nonumber \\&&
\quad =
\sum_{(n_1,\dots,n_l) \in \mathbb L_l}
{(-1)^{\sum_{r=1}^l n_r} \over n_1! \cdots n_l !}
{\partial^{\sum_{j=1}^l n_j} \over
\partial \xi^{n_1} \cdots \partial \xi^{n_l} }
\lim_{Q \to \infty}
\Big \langle \prod_{j=1}^\infty
\Big (1 - \sum_{r=1}^k \xi_r \chi_{\lambda_j - j \in I_r}
\Big ) \Big \rangle^{(\alpha, Q)} 
\Big |_{\alpha = 1 - 2 w/Q^{1/6}}. \nonumber \\
\end{eqnarray}
The limit in the final expression is evaluated according to
(\ref{kt.2}), which shows (\ref{f1f}) can be written
\begin{equation}\label{f2f}
\sum_{(n_1,\dots,n_l) \in \mathbb L_l}
E^{(w)}(\{n_r,(s_r,s_{r-1})\}_{r=1,\dots,l}),
\end{equation}
where $E^{(w)}(\{n_r,(s_r,s_{r-1})\}_{r=1,\dots,l})$ denotes the probability
that the interval $(s_r,s_{r-1})$ contains no particles in the point
process specified by the correlations (\ref{1.12}) with $w$ given by 
(\ref{1.10a}). But (\ref{f2f}) is, according to (\ref{6a.c}), equal to the
distribution function $F^{\symmO}(s_1,\dots,s_l;w)$ 
appearing on the right hand side of
(\ref{1.11}).

\section*{Acknowledgement}

One of the authors (T.N.) is grateful to Dr.~Tomohiro Sasamoto 
for valuable discussions. The referee is to be thanked for some
useful suggestions. The work of PJF was supported by the
Australian Research Council.
 
\bibliographystyle{plain}

\end{document}